\documentclass[a4paper, 11pt, reqno]{amsart}
\usepackage{a4wide}
\usepackage{amssymb, amsmath, amsthm}

\usepackage{hyperref}

\allowdisplaybreaks

\newtheorem{de}{Definition}[section]
\newtheorem{lem}[de]{Lemma}
\newtheorem{prop}[de]{Proposition}

\newtheorem{thm}[de]{Theorem}

\newtheorem{rem}[de]{Remark}

\numberwithin{equation}{section}

\newcommand{\BlackBox}{\rule{1.5ex}{1.5ex}}
\newenvironment{prf}{\par\noindent{\bf Proof.\ }}{\hfill\rule{1.5ex}{1.5ex}\vspace{0.3cm}}

\newcommand{\bprf}{\begin{prf}}
\newcommand{\eprf}{\end{prf}}

\usepackage{amssymb}
\usepackage{dsfont}
\usepackage[T1]{fontenc}

\newcommand{\R}{\mathbb{R}}

\newcommand{\C}{\mathbb{C}}
\newcommand{\N}{\mathbb{N}}

\newcommand{\E}{\mathbb{E}}
\newcommand{\F}{\mathbb{F}}

\newcommand{\X}{\mathbb{X}}
\newcommand{\Y}{\mathbb{Y}}

\newcommand{\D}{\,\textrm{d}}

\newcommand{\trm}{\textrm}

\newcommand{\bea}{\begin{eqnarray*}}
\newcommand{\eea}{\end{eqnarray*}}

\newcommand{\beq}{\begin{equation}}
\newcommand{\eeq}{\end{equation}}

\newcommand{\1}{\mathds{1}}
\newcommand{\tr}{\text{tr}}

\newcommand{\ra}{\rightarrow}

\newcommand{\hra}{\hookrightarrow}

\newcommand{\bc}{\begin{mathcal}}
\newcommand{\ec}{\end{mathcal}}

\newcommand{\calE}{\bc E\ec}

\newcommand{\calH}{\bc H\ec}

\newcommand{\calT}{\bc T\ec}
\newcommand{\calB}{\bc B\ec}
\newcommand{\calD}{\bc D\ec}
\newcommand{\calR}{\bc R\ec}

\newcommand{\HT}{\calH\calT}

\newcommand{\ol}{\overline}
\newcommand{\vphi}{\varphi}

\newcommand{\wt}{\widetilde}
\newcommand{\bs}{\backslash}

\newcommand{\eps}{\varepsilon}

\def\typeout#1{\message{^^J}\message{#1}\message{^^J}}
\newif\ifSRCOK \SRCOKtrue
\newcount\PAGETOP
\newcount\LASTLINE
\global\PAGETOP=1
\global\LASTLINE=-1
\def\EJECT{\SRC\eject}
\def\WinEdt#1{\typeout{:#1}}
\gdef\MainFile{\jobname.tex}
\gdef\CurrentInput{\MainFile}
\newcount\INPSP
\global\INPSP=0
\def\SRC{\ifSRCOK%
  \ifnum\inputlineno>\LASTLINE%
    \ifnum\LASTLINE<0%
      \global\PAGETOP=\inputlineno%
    \fi%
    \global\LASTLINE=\inputlineno%
    \ifnum\INPSP=0%
      \ifnum\inputlineno>\PAGETOP%
        
      \fi%
    \else%
      
    \fi%
  \fi%
\fi}
\def\PUSH#1{%
\SRC%
\ifnum\INPSP=0 \global\let\INPSTACKA=\CurrentInput \else%
\ifnum\INPSP=1 \global\let\INPSTACKB=\CurrentInput \else%
\ifnum\INPSP=2 \global\let\INPSTACKC=\CurrentInput \else%
\ifnum\INPSP=3 \global\let\INPSTACKD=\CurrentInput \else%
\ifnum\INPSP=4 \global\let\INPSTACKE=\CurrentInput \else%
\ifnum\INPSP=5 \global\let\INPSTACKF=\CurrentInput \else%
               \global\let\INPSTACKX=\CurrentInput \fi\fi\fi\fi\fi\fi%
\gdef\CurrentInput{#1}%
\WinEdt{<+ \CurrentInput}%
\global\LASTLINE=0%
\ifSRCOK\fi%
\global\advance\INPSP by 1}
\def\POP{%
\ifnum\INPSP>0 \global\advance\INPSP by -1  \fi%
\ifnum\INPSP=0 \global\let\CurrentInput=\INPSTACKA \else%
\ifnum\INPSP=1 \global\let\CurrentInput=\INPSTACKB \else%
\ifnum\INPSP=2 \global\let\CurrentInput=\INPSTACKC \else%
\ifnum\INPSP=3 \global\let\CurrentInput=\INPSTACKD \else%
\ifnum\INPSP=4 \global\let\CurrentInput=\INPSTACKE \else%
\ifnum\INPSP=5 \global\let\CurrentInput=\INPSTACKF \else%
               \global\let\CurrentInput=\INPSTACKX \fi\fi\fi\fi\fi\fi%
\WinEdt{<-}%
\global\LASTLINE=\inputlineno%
\global\advance\LASTLINE by -1%
\SRC}
\def\INPUT#1{\relax}

\let\originalxxxeverypar\everypar
\newtoks\everypar
\originalxxxeverypar{\the\everypar\expandafter\SRC}
\everymath\expandafter{\the\everymath\expandafter\SRC}
\output\expandafter{\expandafter\SRCOKfalse\the\output}

\newif\ifSRCOK \SRCOKtrue
\newcount\PAGETOP
\newcount\LASTLINE
\global\PAGETOP=1
\global\LASTLINE=-1
\gdef\MainFile{\jobname.tex}
\gdef\CurrentInput{\MainFile}
\newcount\INPSP
\global\INPSP=0
\def\EJECT{\SRC\eject}
\def\WinEdt#1{\typeout{:#1}}
\def\SRC{\ifSRCOK%
  \ifnum\inputlineno>\LASTLINE%
    \ifnum\LASTLINE<0%
      \global\PAGETOP=\inputlineno%
    \fi%
    \global\LASTLINE=\inputlineno%
    \ifnum\INPSP=0%
      \ifnum\inputlineno>\PAGETOP%
      \fi%
    \else%
    \fi%
  \fi%
\fi}
\def\PUSH#1{%
\SRC%
\ifnum\INPSP=0 \global\let\INPSTACKA=\CurrentInput \else%
\ifnum\INPSP=1 \global\let\INPSTACKB=\CurrentInput \else%
\ifnum\INPSP=2 \global\let\INPSTACKC=\CurrentInput \else%
\ifnum\INPSP=3 \global\let\INPSTACKD=\CurrentInput \else%
\ifnum\INPSP=4 \global\let\INPSTACKE=\CurrentInput \else%
\ifnum\INPSP=5 \global\let\INPSTACKF=\CurrentInput \else%
               \global\let\INPSTACKX=\CurrentInput \fi\fi\fi\fi\fi\fi%
\gdef\CurrentInput{#1}%
\WinEdt{<+ \CurrentInput}%
\global\LASTLINE=0%
\ifSRCOK\fi%
\global\advance\INPSP by 1}
\def\POP{%
\ifnum\INPSP>0 \global\advance\INPSP by -1  \fi%
\ifnum\INPSP=0 \global\let\CurrentInput=\INPSTACKA \else%
\ifnum\INPSP=1 \global\let\CurrentInput=\INPSTACKB \else%
\ifnum\INPSP=2 \global\let\CurrentInput=\INPSTACKC \else%
\ifnum\INPSP=3 \global\let\CurrentInput=\INPSTACKD \else%
\ifnum\INPSP=4 \global\let\CurrentInput=\INPSTACKE \else%
\ifnum\INPSP=5 \global\let\CurrentInput=\INPSTACKF \else%
               \global\let\CurrentInput=\INPSTACKX \fi\fi\fi\fi\fi\fi%
\WinEdt{<-}%
\global\LASTLINE=\inputlineno%
\global\advance\LASTLINE by -1%
\SRC}
\def\INPUT#1{\relax}
\let\OldINCLUDE=\include
\def\include#1{
\EJECT%
\PUSH{#1.tex}%
\OldINCLUDE{#1}%
\POP}

\let\originalxxxeverypar\everypar
\newtoks\everypar
\originalxxxeverypar{\the\everypar\expandafter\SRC}
\everymath\expandafter{\the\everymath\expandafter\SRC}
\let\zzzxxxbibliography=\bibliography
\def\bibliography#1{\PUSH{\jobname.bbl}\zzzxxxbibliography{#1}\POP}
\output\expandafter{\expandafter\SRCOKfalse\the\output}

\begin{document}

\title{Interpolation, embeddings and traces of anisotropic fractional
Sobolev spaces with temporal weights}

\author{Martin Meyries}
\author{Roland Schnaubelt}

\address{Department of Mathematics,
Karlsruhe Institute of Technology, 76128 Karlsruhe, Germany.}
\email{martin.meyries@kit.edu}
\email{schnaubelt@kit.edu}

\thanks{This paper is part of a research project supported by the Deutsche
Forschungsgemeinschaft (DFG)}




\begin{abstract} We investigate  the properties of a class of weighted vector-valued $L_p$-spaces and the corresponding (an)isotropic Sobolev-Slobodetskii spaces. These spaces arise naturally in the context of maximal $L_p$-regularity for parabolic initial-boundary value problems. Our  main tools are operators with a bounded $\calH^\infty$-calculus, interpolation theory, and operator sums. \end{abstract}

\maketitle

\section{Introduction and preliminaries}

In this paper we investigate a class of anisotropic fractional Sobolev spaces on space-time
with weights in the time variable. We treat in particular interpolation results for such spaces, 
Sobolev type embeddings as well as temporal and spatial trace theorems
in a systematic, operator theoretic way. Our choice of 
both the class of spaces and the studied properties is motivated by applications to 
quasilinear parabolic evolution equations with nonlinear boundary conditions which will be treated
in subsequent papers focussing on the longterm behavior, cf.\ \cite{Mey10}. These papers will be 
based on linearization; i.e.,
one writes the nonlinear equation as a linear initial-boundary value problem whose inhomogeneities
contain the nonlinear terms.  Since the underlying problem is quasilinear, the linear and the 
nonlinear parts will be of the same order so that it is crucial to have sharp results for the
linear part giving optimal regularity. In the companion paper \cite{MS11} we establish the necessary
 theory for linear inhomogeneous initial-boundary value problems,
based on the present study of the underlying function spaces.

The main focus of this work is the presence of the temporal weights. To explain their role,
we first recall some features of the known theory for the unweigthed case by means of an 
example. For a domain $\Omega\subset \R^n$ with smooth compact boundary $\partial \Omega$
and given functions 
$f$, $g$ and $u_0$, we consider the heat equation
\begin{alignat}{3} \label{heat}
\partial_t u(t,x) - \Delta u(t,x) & =   f(t,x), &\qquad &  x\in\Omega,   &\qquad & t>0,  \notag \\
\partial_\nu u(t,x)   & =  g(t,x),  && x\in\partial \Omega, &&  t>0, \\
u(0,x) & = u_0(x), &&  x\in\Omega, && \notag
\end{alignat}
with an inhomogeneous Neumann boundary condition. We want to work in a framework where
the boundary conditions can be understood in the sense of traces, and not just weakly.
So we choose an $L_p$ setting and require that the solution $u$ belongs to the space
$$\E_{1}:= W_{p}^1\big (J; L_p(\Omega)\big)\cap L_{p}\big(J;W^2_p(\Omega)\big)$$ 
for some $p\in(1,\infty)$ and  a finite interval $J=(0,T)$. 
To obtain the optimal spaces for the data, one needs sharp trace 
results for the space $\E_1$. It is known that the Neumann trace $g$ of $u\in\E_1$ belongs to 
the anisotropic Slobodetskii space
$$\F_1=W_{p}^{1/2-1/2p}\big(J; L_p(\partial\Omega)\big)\cap 
        L_p\big(J; W_p^{1-1/p}(\partial\Omega)\big)$$ 
and that the trace $u_0$ of $u$ at time $t=0$ belongs to $W^{2-2/p}_p(\Omega)$. Observe that 
$\F_1$ still retains some time regularity! Moreover, one has $f\in L_p\big (J; L_p(\Omega)\big )$ 
and the compatibility condition $\partial_\nu u_0 =g|_{t=0}$
should hold at time $t=0$ if $g$ has a trace at $t=0$, which happens for $p>3$. In fact,
 problem \eqref{heat} has a unique solution $u\in\E_1$ if and only if
these conditions hold. We refer to e.g.\ \cite{DHP07}
 for these and much more general results. Further,  the trace space $W^{2-2/p}_p(\Omega)$ 
 is the smallest space  in which all solutions  
 of \eqref{heat} in $\E_1$ are continuous on $[0,T]$, and the corresponding nonlinear problems share
  this property. Correspondingly, the (local) semiflows of solutions live
in $W^{2-2/p}_p(\Omega)$ and  the norm of this space is the right one to describe
their properties. For instance, it gives the  blow-up condition, see e.g.\ \cite{LPS06}.
One often takes a large $p$ to simplify the treatment of the nonlinearities (say, $p>n+2$).
 Thus the norm of $W^{2-2/p}_p(\Omega)$ is far away from the norms one can control by 
 the usual a priori estimates, such as the norms of $W^1_2$ or  $L_\infty$. This unpleasant
regularity gap is even more significant in the  H\"older setting.

There is a known  technique to reduce the necessary  regularity for initial
conditions of evolution equations. One studies the problem in the weighted spaces
$$L_{p,\mu}(\R_+;E) := \big \{u:\R_+\ra E \,:\, t^{1-\mu}u\in  L_p(\R_+;E)\big\},$$
where $E$ is a Banach space  and  $\mu\in(1/p,1]$.
The corresponding weighted Sobolev spaces are defined similarly, and 
the fractional Sobolev spaces $W_{p,\mu}^s(\R_+;E)$ and $H^s_{p,\mu}(\R_+;E)$ are introduced
via real and complex interpolation, respectively. 
Observe that for $\mu=1$ we recover the unweigthed case and that the spaces become larger
if we decrease $\mu$. So this type of weights allows for functions being more singular at $t=0$. 
(We note that the weight $t^{p(1-\mu)}$ belongs to the class $A_p$ used  in harmonic analysis, 
see e.g.\ \cite{St93}.) 

This approach has been carried out in an $L_p$ setting by Pr\"uss $\&$ Simonett
in \cite{PS04} for problems with homogeneous boundary conditions (i.e., $g=0$ in \eqref{heat}). 
In this case one can work entirely in the framework of semigroup theory. In fact, for the generator $-A$
of an exponentially stable semigroup one looks for a solution $u$ in 
$$\E_{1,\mu}:= W_{p,\mu}^1(\R_+; E)\cap L_{p,\mu}(\R_+;D(A))$$ 
for the evolution equation
\begin{equation}\label{homog} 
 u'(t) + Au(t) = f(t), \quad t>0, \qquad u(0) = u_0.
\end{equation}
Clearly, it must hold that $f\in L_{p,\mu}(\R_+;E)$. It is known 
that the initial value $u_0$ has to belong to the real interpolation space 
$(E, D(A))_{\mu-1/p,p}$,  see e.g.\ Theorem~1.14.5 in \cite{Tri94}.
The main result in \cite{PS04} now says that one has a unique solution  $u\in\E_{1,\mu}$ of 
\eqref{homog} for all $f\in L_{p,\mu}(\R_+;E)$ and $u_0\in (E, D(A))_{\mu-1/p,p}$
if (and only if) $A$ has maximal $L_p$-regularity; i.e., \eqref{homog} with $u_0=0$ has a solution
$u\in\E_1$ for each $f\in   L_{p}(\R_+;E)$. The latter property is well understood, see 
\cite{KW04} and the references therein. As a result, one can reduce the needed initial 
regularity for \eqref{homog} almost to the base space, say, $E=L_p(\Omega)$.
We further point out  that the temporal weights give 
the scale of phase spaces  $(E, D(A))_{\mu-1/p,p}$, which are compactly embedded into each other
in many cases. Precisely this point was used  in the very recent paper \cite{KPW10}
to establish attractivity of equilibria of quasilinear evolution equations based on the
results of \cite{PS04}, see also \cite{KPW11}.

Unfortunately, it seems to be impossible obtain sharp (Sobolev) regularity for inhomogeneous 
boundary value problems such as \eqref{heat} in a pure semigroup framework. Instead one
has to restrict to a PDE setting and reduce equations like \eqref{heat} to model problems on full- 
and half-spaces by means of local charts  and the freezing of coefficients. Nevertheless
 we see below that semigroup and operator theory can still play a crucial 
role in the proofs. In this way, vector-valued parabolic initial-boundary value problems with 
inhomogeneous boundary conditions were treated by Denk, Hieber, Pr\"uss $\&$ Zacher in the papers 
\cite{DHP03}, \cite{DHP07} and \cite{DPZ08}  without temporal weights. In the present and 
the following paper \cite{MS11} we extend their results to the case of the weights $t^{p(1-\mu)}$, 
focussing here on the function spaces themselves. Our setting also covers linearizations of 
problems arising from free boundary problems, see e.g.\ \cite{EPS03} and \cite{KPW11}.

In Section~\ref{sec:basic} we start with basic properties of the weighted spaces.
The crucial point is that the operator $-\partial_t$ is maximally accretive and  
has a bounded $\calH^\infty$-calculus on $L_{p,\mu}(\R_+;E)$ (if $E$ is of class $\calH\calT$), 
see Lemma~\ref{sec:timeder} and  Theorem~\ref{wm11}. Using a theorem by Yagi,  one can 
conclude that the spaces $H^{s}_{p,\mu}(\R_+;E)$ are the domains of fractional powers of 
$1-\partial_t$, see 
Proposition~\ref{maal1}. This fact also allows us to establish the natural interpolation 
properties for the scales of  $W^{s}_{p,\mu}$ and $H^{s}_{p,\mu}$ in Lemma~\ref{lg54}. 
These results are complemented by several propositions on extension, density,  intrinsic
norms, Sobolev embeddings 
and a Poincar\'e inequality. The proofs rely on interpolation theory, Hardy's inequality and 
results by Grisvard, \cite{Gri63}. 

In view of the above considerations it is clear that we have to study the mapping properties
of temporal and spatial derivatives and traces  in anisotropic spaces like $\E_{1,\mu}$. 
Here spaces like
$$W_{p,\mu}^{s+\alpha}\big (\R_+; W_p^r(\Omega;E)\big)\cap 
  W_{p,\mu}^s\big(\R_+; W_p^{r+\beta}(\Omega;E)\big), \qquad s,r,\alpha,\beta \geq 0,$$ 
naturally appear, as well as analogous Bessel potential spaces $H$.  Of course, one also needs
Sobolev type embeddings. By localization and local charts, we reduce these results to the model
cases of full- and half-spaces. There we use that the related $H$ spaces are the domains
of the operator $L=(1-\partial_t)^\alpha+(1-\Delta)^{\beta/2} $ having
 bounded imaginary powers, see Lemma~\ref{wm22}. Then Sobolevskii's
mixed derivative theorem and interpolation arguments imply the fundamental embedding result 
Proposition~\ref{embani} which allows to interpolate between time and space regularity. 
We can conclude the desired mapping properties of spatial derivatives in Lemma~\ref{spatiald}.
Another crucial ingredient is Lemma~\ref{lg72} on time traces of semigroup orbits in weighted
spaces. Employing these tools, 
we then establish our main Theorem~\ref{mott} on the time traces of the  anisotropic fractional
spaces. Similarly, we deduce the final  Theorem~\ref{sec:generalstt} on spatial traces.
It turns out that our main results are natural extensions of the unweighted case $\mu=1$, and 
thus there is no disadvantage when working in a weighted framework in the context of 
parabolic problems.

The basic strategy of our proofs of the trace theorems is similar to that in e.g.\ 
Section~3 of \cite{DHP07} where the unweighted case was treated, see also \cite{EPS03} or 
\cite{Zac03}. The more recent paper \cite{Ama09} gives a comprehensive account 
of unweighted anisotropic spaces. However, for the weighted case we first had to establish 
the underlying theory contained in Sections~\ref{sec:basic} and \ref{was} and in Lemma~\ref{lg72}. 
As a fundamental difference to the unweighted case, our weighted spaces are not invariant under 
the right shift and have no suitable extension to the full time interval $\R$. This behavior 
reflects the desired fact that our weights vanish at $t=0$ and allow for stronger singularities 
at $t=0$. Moreover, some proofs in the literature cannot simply be generalized to the weighted case 
(e.g., those of the temporal trace theorems in \cite{EPS03} or \cite{Zac03} which use fractional 
evolution equations).

We also point out that the multiplication operator $\Phi_\mu u= t^{1-\mu} u$ does \textsl{not}
give an isomorphism from $W^1_{p,\mu}(\R_+;E)$ to $W^1_{p}(\R_+;E)$, since 
$(\Phi_\mu\1)'\notin L_p(0,1)$. As a consequence, one cannot simply reduce the results to the 
unweighted case by isomorphy. For other types of weights this approach would be possible, 
for instance for exponential weights $e^{-p\rho t}$ which do not change the behavior at $t=0$.
This was done in the recent paper \cite{DSS08} where interpolation and trace theorems are
shown for the exponentially weighted spaces. (Of course, this paper focusses on
different results.) Isotropic weighted spaces have been studied in detail in the literature,
see e.g.\ \cite{Gri63}, \cite{Kuf85} and \cite{Tri94}. 

In this work we strive for a self contained and systematic treatment of the results needed
for evolution equations. This should be of some interest also in the known unweighted case 
$\mu=1$. To keep the length within reasonable bounds, we omitted some arguments or calculations
which are either routine or similar to parallel cases discussed here or in the literature. 
Most of the omitted details are presented in the Ph.D.\ thesis \cite{Mey10}
to which we then refer.

\medskip

We now indicate some concepts and results from operator theory that we
use frequently. Details can be found at the references  given below.
 The real and complex interpolation functors  $(\cdot,\cdot)_{\theta,q}$
and $[\cdot,\cdot]_\theta$, repectively, are a fundamental tool in our investigations. 
For the relevant theory the reader is referred to \cite{Lun09} and \cite{Tri94}.
If $A$ is a sectorial operator on a Banach space $E$, $\theta\in (0,1)$ and $q\in [1,\infty]$, 
 we set  $D_A(\theta,p) := \big( E,D(A)\big)_{\theta,q}$.
The theory of strongly continuous semigroups is presented in e.g.\ \cite{EN00}.
Throughout we assume that the Banach space $E$ is of class $\calH\calT$ 
(or, equivalently, is a UMD space), 
which means that the Hilbert transform is bounded on $L_2(\R;E)$, see Sections~III.4.3--5 
of \cite{Ama95}. We note that Hilbert spaces are   of class $\calH\calT$, as well
as the reflexive Lebesgue and (fractional)  Sobolev spaces.
The definition and properties of the bounded
$\calH^\infty$-calculus of a sectorial operator $A$ are discussed in  \cite{DHP03} or \cite{KW04}.
However, in the sequel we almost never use the calculus itself, but rather
some theorems requiring it as an assumption. An important consequence of the 
 bounded $\calH^\infty$-calculus is that $A$ admits bounded imaginary powers,
see e.g.\ Sections~2.3 and 2.4 of \cite{DHP03} or Section~III.4.7 of \cite{Ama95}. 
Finally, we recall the Dore-Venni theorem from \cite{DV87} in a version of \cite{PrS90} and 
combined with the mixed derivative theorem due to Sobolevskii, \cite{Sob75}.
\begin{prop}\label{dv}
\textsl{Let $E$ be a Banach space of class $\calH\calT$ and suppose that the operators $A,B$ 
on $E$ are resolvent commuting and admit bounded imaginary powers with power angles satisfying 
$\theta_A+\theta_B < \pi$. Assume further that $A$ or $B$ is invertible. Then $A+B$ is invertible,
 admits bounded imaginary powers with angle not larger than $\theta_A+\theta_B$ and  the operator 
 $A^\alpha B^{1-\alpha} (A+B)^{-1}$ is bounded on $E$ for all $\alpha\in [0,1]$.}
 \end{prop}We write $a\lesssim b$ for some quantities $a,b$ if there is a 
generic positive constant $C$ with $a\leq Cb$. If $X,Y$ are Banach spaces, $\calB(X,Y)$ 
is the space of bounded linear operators between them, and we put $\calB(X):= \calB(X,X)$.

\section{Basic properties of the weighted spaces}\label{sec:basic}
Throughout, we consider a time interval $J= \R_+:= (0,\infty)$ or $J=(0,T)$ for some $T>0$
and a Banach space $E$ of class $\calH\calT$, and we let 
$$p\in (1,\infty), \qquad \mu\in (1/p,1].$$ 
We study the Banach space $L_{p,\mu}(J;E) := \{u:J\ra E \,:\, t^{1-\mu}u\in  L_p(J;E) \}$
with the norm
$$|u|_{L_{p,\mu}(J;E)} := |t^{1-\mu} u|_{L_p(J;E)} 
  = \left ( \int_J t^{p(1-\mu)}|u(t)|_E^p \D t\right )^{1/p}.$$
 Occasionally, we use these spaces also with $\mu$ replaced by some $\alpha>1/p$.
It clearly holds
$$L_p(0,T;E)\hra L_{p,\mu}(0,T;E) \quad \text{ and } \quad
 L_{p,\mu}(0,T;E)\hra L_p(\tau,T;E)$$
 for all $\tau\in (0,T)$, but $L_p(\R_+;E) \nsubseteq L_{p,\mu}(\R_+;E)$ for 
 $\mu\in (1/p,1)$. For $k\in \N_0$, we define the corresponding weighted Sobolev space
$$ W_{p,\mu}^k(J;E) = H_{p,\mu}^k(J;E) := \big \{u\in W_{1,\text{loc}}^k(J;E) \;:\; u^{(j)} \in  
  L_{p,\mu}(J;E),\;\; j\in \{0,...,k\} \big \}$$ 
(where $ W_{p,\mu}^0=H_{p,\mu}^0=L_{p,\mu}$ by definition), 
which is a Banach space endowed with the norm 
$$ |u|_{W_{p,\mu}^k(J;E)} = |u|_{H_{p,\mu}^k(J;E)} 
     :=  \bigg ( \sum_{j=0}^k |u^{(j)}|_{L_{p,\mu}(J;E)}^p \bigg )^{1/p}.$$ 
For $s\in \R_+\bs\N$,  we write $s= [s] + s_*$ with $[s]\in \N_0$ and $s_*\in (0,1)$.
The weighted Slobodetskii and Bessel potential spaces 
$$W_{p,\mu}^s(J;E)  :=   \big ( W_{p,\mu}^{[s]}(J;E), W_{p,\mu}^{[s]+1}(J;E)  \big)_{s_*,p}\,,\quad
H_{p,\mu}^s(J;E)  :=   \big[W_{p,\mu}^{[s]}(J;E), W_{p,\mu}^{[s]+1}(J;E)  \big]_{s_*}$$
are introduced by means of real and complex interpolation, respectively.
In view of Theorem~4.2/2 of \cite{Sch} and Satz~3.21 of \cite{Zi89}, this definition is 
consistent with the unweighted case; i.e., we have $W_p^s = W_{p,1}^s$ and $H_p^s = H_{p,1}^s$ 
for all $s\geq 0$. The general properties of real and complex interpolation spaces 
(see \cite{Lun09} or  \cite{Tri94}) imply that  one has the scale of dense embeddings
\beq\label{kv33}
W_{p,\mu}^{s_1}\stackrel{d}{\hra} H_{p,\mu}^{s_2} \stackrel{d}{\hra} 
W_{p,\mu}^{s_3}\stackrel{d}{\hra} H_{p,\mu}^{s_4}, \qquad s_1>s_2>s_3>s_4\geq 0,
\eeq
for all fixed $p\in (1,\infty)$ and $\mu\in (1/p,1]$. In the sequel we will often use that 
$$\calB\big (W_{p,\mu}^{k}(J;E)\big ) \cap \calB\big(W_{p,\mu}^{k+1}(J;E)\big) \hra 
\calB\big(W_{p,\mu}^{s}(J;E)\big)\cap \calB\big(H_{p,\mu}^{s}(J;E)\big),$$
where $k\in \N_0$ and $s\in (k,k+1)$, to deduce assertions for exponents $s\ge0$ 
from the integer case. Before giving further definitions, we state a first basic property 
which can easily be proved using H\"older's inequality.

\begin{lem}\label{wm14}
\textsl{Let $J=(0,T)$ be finite or infinite, $p\in (1,\infty)$ and $\mu \in (1/p,1]$. 
We then have 
\begin{equation*}
 L_{p,\mu}(J;E) \hra L_{q,\emph{\text{loc}}}(\ol{J};E)
 \qquad \text{for all \ } 1\leq q< \frac{1}{1-\mu+1/p}\le p.\tag*{\BlackBox}
\end{equation*}}
\end{lem}

\noindent
It thus holds $W_{p,\mu}^k(J;E) \hra W_{1,\text{loc}}^k(\ol{J};E)$
for all $k\in \N$ so that the trace $u\mapsto u^{(j)}(0)$ is continuous
from $W_{p,\mu}^k(J;E)$ to $E$ for all $j\in \{0,...,k-1\}$. So we can define 
$${}_0W_{p,\mu}^k(J;E) = {}_0H_{p,\mu}^k(J;E) := \big \{ u\in W_{p,\mu}^k(J;E)\;:\; 
u^{(j)}(0) = 0,\;\; j\in \{0,...,k-1\}\big \}$$ 
for $k\in \N$ which are Banach spaces with the norm of $W_{p,\mu}^k$. 
For convenience we further set 
$${}_0W_{p,\mu}^0(J;E) = {}_0H_{p,\mu}^0(J;E) := L_{p,\mu}(J;E).$$
For a number $s=[s] +s_*\in \R_+\bs \N$, we again define the corresponding fractional 
order spaces by interpolation; i.e., we put
\begin{align*}
{}_0W_{p,\mu}^s(J;E) & :=   \big({}_0W_{p,\mu}^{[s]}(J;E), {}_0W_{p,\mu}^{[s]+1}(J;E)  
                            \big)_{s_*,p}\,,\\
{}_0H_{p,\mu}^s(J;E) & :=   \big[{}_0W_{p,\mu}^{[s]}(J;E), {}_0W_{p,\mu}^{[s]+1}(J;E)  
                           \big]_{s_*}\,,
\end{align*}                           
obtaining as before a scale of function spaces 
\beq\label{kv5000}
{}_0 W_{p,\mu}^{s_1} \stackrel{d}{\hra} {}_0 H_{p,\mu}^{s_2} \stackrel{d}{\hra} 
   {}_0 W_{p,\mu}^{s_3}\stackrel{d}{\hra} {}_0 H_{p,\mu}^{s_4}, \qquad s_1>s_2>s_3>s_4\geq 0.
\eeq 
We further have
$${}_0W_{p,\mu}^s(J;E) \hra W_{p,\mu}^s(J;E) \quad \text{ and  } \quad
{}_0H_{p,\mu}^s(J;E) \hra H_{p,\mu}^s(J;E) \qquad \text{for all \ } s\ge 0.$$
 Hardy's inequality  is a fundamental tool for the 
 investigation of the $L_{p,\mu}$-spaces. It states that
\beq\label{hardy-ineq}
\int_0^\infty \left( t^{-\alpha} \int_0^t \varphi(\tau)\D \tau\right)^p \D t 
 \leq  \frac{1}{(\alpha-1/p)^p} \int_0^\infty (t^{1-\alpha} \varphi(t))^p \D t
 \eeq 
for all $\alpha \in (1/p,\infty)$  and measurable nonnegative functions $\varphi:\R_+\to\R$,
see Theorem~330 in \cite{HLP34}.
For the spaces with vanishing initial values we deduce the following weighted inequality.

\begin{lem} \label{sec:Hardy}
\textsl{Let $J=(0,T)$ be finite or infinite, $p\in (1,\infty)$, $\mu \in (1/p,1]$ 
and $s\geq 0$. Then 
$$\int_J t^{p(1-\mu - s)} |u(t)|_E^p \D t \leq C_{p,\mu,s} \,|u|_{\,{}_0W_{p,\mu}^s(J;E)}^p
 \qquad \text{if }\;u\in {}_0W_{p,\mu}^s(J;E),$$ 
and this remains true if one replaces  ${}_0W_{p,\mu}^s(J;E)$ by ${}_0H_{p,\mu}^s(J;E)$.}
\end{lem}
\bprf The case $s=0$ is trivial and the case $s=1$ is
a consequence of Hardy's inequality with $\varphi= |u'|_E$.  The assertion for $s=k\in\N_0$
follows by induction, so that we have 
\beq\label{lg51}
{}_0W_{p,\mu}^k(J;E) \hra L_p(J, t^{p(1-\mu-k)}\D t; E).
\eeq 
Theorem~1.18.5 of \cite{Tri94} yields that
\beq\label{lg50}
\big (L_p(J,t^{p(1-\mu-k)}\D t; E), L_p(J,t^{p(1-\mu-(k+1))}\D t; E)\big )_{\theta,p}  
= L_p(J,t^{p(1-\mu-\theta k)}\D t; E)
\eeq
for all $\theta\in (0,1)$ and that (\ref{lg50}) remains true if one replaces
 $(\cdot,\cdot)_{\theta,p}$ by $[\cdot,\cdot]_\theta$. Interpolation of the embedding (\ref{lg51})
 now implies  the assertions for all $s\ge0$.
\eprf

The above inequalities can be used to show that the multiplication 
operator $\Phi_\mu u= t^{1-\mu} u$ is an isomorphism from the weighted to the
unweighted spaces, provided one restricts to vanishing initial values. Without this 
restriction, $\Phi_\mu$ does not map properly since $(\Phi_\mu\1)'\notin L_p$.

\begin{lem}\label{isophi}
\textsl{Let $J=(0,T)$ be finite or infinite, $p\in (1,\infty)$, $\mu\in (1/p,1]$ 
and $s\geq 0$. Then the map $\Phi_\mu$, given by 
$$(\Phi_\mu u)(t) := t^{1-\mu}u(t),$$ 
induces an isomorphism from ${}_0W_{p,\mu}^s(J;E)$ to ${}_0W_{p}^s(J;E)$ and from
 ${}_0H_{p,\mu}^s(J;E)$ to ${}_0H_{p}^s(J;E)$. The inverse $\Phi_\mu^{-1}$ 
 is given by $(\Phi_\mu^{-1} u)(t) = t^{-(1-\mu)} u(t)$.\hfill\rule{1.5ex}{1.5ex}}
\end{lem}
The lemma is shown for $s\in \{0,1\}$ in Proposition~2.2 of \cite{PS04}. The case of other integers
can be treated similarly, see Lemma~1.1.3 of \cite{Mey10}. The general case then follows by 
interpolation. We use the above lemma to prove basic density results for the weighted spaces.

\begin{lem}\label{wm7}
\textsl{Let $J=(0,T)$ be finite or infinite, $p\in (1,\infty)$, $\mu \in (1/p,1]$ and $s\geq 0$. 
Then the space $C_c^\infty(\ol{J}\bs\{0\};E)$ is dense in ${}_0W_{p,\mu}^s(J;E)$
and  ${}_0H_{p,\mu}^s(J;E)$, and $C_c^\infty(\ol{J};E)$ is dense in $W_{p,\mu}^s(J;E)$ 
and $H_{p,\mu}^s(J;E).$}
\end{lem}
\bprf Due to (\ref{kv33}) and (\ref{kv5000}) we only have to consider the case $s=k\in \N_0$. 
The preceding lemma  shows that $\Phi_\mu u$ belongs to  ${}_0W_{p}^k(J;E)$ for every
 $u\in {}_0W_{p,\mu}^k(J;E)$. As in Theorem 2.9.1 of \cite{Tri94} for the scalar-valued case, 
one sees that $C_c^\infty(\ol{J}\bs\{0\};E)$ is dense in ${}_0W_{p}^k(J;E)$ for $k\in \N_0$. 
For any given $\eps>0$ there thus exists a function $\psi\in C_c^\infty(\ol{J}\bs\{0\};E)$ 
with $|\Phi_\mu u - \psi|_{W_{p}^k(J;E)} < \eps$. Consequently,
$|u - \Phi_\mu^{-1} \psi|_{W_{p,\mu}^k(J;E)} \lesssim \eps$,
and so $C_c^\infty(\ol{J}\bs\{0\};E)$ is dense in ${}_0W_{p,\mu}^k(J;E)$. 

For $u\in W_{p,\mu}^k(J;E)$, we choose $\psi_1\in C_c^\infty(\ol{J};E)$ with 
$\psi_1^{(j)}(0) = u^{(j)}(0)$ for $j\in \{0,...,k-1\}$. Since $u-\psi_1\in {}_0W_{p,\mu}^k(J;E)$, 
by the above considerations
there is $\psi_2\in C_c^\infty(\ol{J}\setminus\{0\};E)$ 
being $\eps$-close to $u-\psi_1$ in ${}_0W_{p,\mu}^k(J;E)$  
Hence, $\psi_1+\psi_2\in C_c^\infty(\ol{J};E)$ is $\eps$-close to $u$ in $W_{p,\mu}^k(J;E)$.
\eprf

We next construct  extension operators to $\R_+$ for the weighted spaces on finite intervals; 
i.e., continuous right-inverses of the restriction operator. For vanishing initial values 
one can achieve an extension whose norm is bounded independently of the length of the 
underlying interval. 
This fact is crucial for the treatment of nonlinear evolution equations.
For simplicity we consider only  $s\in [0,2]$ in this case.

\begin{lem}\label{sec:ext3}
\textsl{Let $J=(0,T)$ be finite and $0<1/p<\mu\le 1$. Given $k\in \N$, there is an extension operator 
$\calE_J$ from $J$ to $\R_+$ with 
$$\calE_J\in \calB\big (W_{p,\mu}^s(J;E), W_{p,\mu}^s(\R_+;E)\big ) \cap \calB \big 
     (H_{p,\mu}^s(J;E), H_{p,\mu}^s(\R_+;E) \big) \quad \text{for all  } s\in[0,k].$$
 Here we can replace $W$ by ${}_0W$ and $H$ by ${}_0H$. There is further an extension operator 
$\calE_J^0$ from $J$ to $\R_+$ with 
$$\calE_J^0\in \calB\big({}_0W_{p,\mu}^s(J;E), {}_0 W_{p,\mu}^s(\R_+;E)\big) \cap 
  \calB({}_0H_{p,\mu}^s(J;E), {}_0H_{p,\mu}^s(\R_+;E)\big) \quad 
\text{for all  }s \in [0,2],$$ whose operator norm is independent of $T$. 
Moreover, $\calE_J, \calE_J^0 \in \calB\big(L_\infty(J;E), L_\infty(\R_+;E)\big)$  
 with operator norms independent of $T$.}
\end{lem}
\bprf 
Let $k\in\N$. We extend $u\in W^k_{p,\mu}(0,T;X)$ to $\calE u
\in W^k_{p,\mu}(0,T_k;X)$ with $T_k=T+T/(2k+2)$ in the same way as in, e.g.,
Theorem~5.19 of \cite{AF03}. In particular, $\calE u(t)$ only depends on $u$ on 
$(T/2,T)$ if $t\ge T$.  Take a function $\varphi\in C^\infty(\R_+)$ which is equal to 1
on $(0,T)$  and has its support in $(0,T_k) $, and set $\calE_J :=\varphi\calE$ (with trivial 
extension for $t\ge T_k$). It is easy to check that $\calE_J$ is bounded from 
$ W^k_{p,\mu}(0,T;X)$ to $W^k_{p,\mu}(\R_+;E)$, and the general case follows from interpolation. 
For $k\in \{0,1,2\}$ and a function $u\in {}_0W_{p,\mu}^k(J;E)$ we define $\calE_J^0$ 
by the formula
$$
(\calE_J^0u)(t) := \begin{cases} u(t), & t\in (0,T), \\ 
    3 (\psi^{1-\mu}  u)(2T-t) \1_{[T, 2T]}(t) - 2 (\psi^{1-\mu}  u)(3T-2t) \1_{[T, 
                     \frac{3}{2}T]}(t), & t\geq T, \end{cases}
$$ 
where $\psi(\tau) = \frac{2T\tau - \tau^2}{T^2}$ and $\1_I$ is the characteristic function 
of an interval $I$. Using Lemma \ref{sec:Hardy}, it can be checked  by a direct calculation 
that $\calE_J^0$ is bounded on ${}_0W_{p,\mu}^k(J;E)$ for all $k\in \{0,1,2\},$
see Lemma~1.1.5 in \cite{Mey10}. Interpolation yields the general case.
 The explicit representations of $\calE_J$ and $\calE_J^0$ show that 
 these operators admit an  $L_\infty$-estimate independent of $T$.\eprf

We now investigate the realization of the derivative $\partial_t$ and its fractional powers 
on the weighted spaces. The properties of these and similar operators are fundamental 
for all our further considerations. We start with a generation result.

\begin{lem}\label{sec:timeder}
\textsl{For $p\in (1,\infty)$, $\mu\in (1/p,1]$ and $s\geq 0$, 
the family of left translations $\{\Lambda_{t}^{E}\}_{t\geq 0}$, given by 
$$(\Lambda_{t}^{E} u)(\tau) := u(\tau+t), \qquad \tau\geq 0,$$ 
forms a strongly continuous contraction semigroup on 
$W_{p,\mu}^s(\R_+;E)$ and on $H_{p,\mu}^s(\R_+;E)$. Its generator is the 
derivative $\partial_t$ with domain $W_{p,\mu}^{s+1}(\R_+; E)$ and $H_{p,\mu}^{s+1}(\R_+; E)$, 
respectively.}
\end{lem}

\bprf\textbf{(I)} We write $\Lambda_t = \Lambda_{t}^{E}$ for simplicity. For each $t_0\geq 0$ the operator $\Lambda_{t_0}$ maps $L_{p,\mu}(\R_+;E)$ into itself and is contractive because of  
$$|\Lambda_{t_0}u|_{L_{p,\mu}(\R_+;E)}^p 
  \leq \int_{0}^\infty (\tau+t_0)^{p(1-\mu)} |u(\tau+t_0)|_E^p\D \tau 
   \leq |u|_{L_{p,\mu}(\R_+;E)}^p.$$
Analoguously one shows that $\Lambda_{t_0}$ is a contractive map on $W_{p,\mu}^k(\R_+;E)$
for all $k\in \N$. By interpolation, this fact carries over to $W_{p,\mu}^s(\R_+;E)$ and
 $H_{p,\mu}^s(\R_+;E)$ for all $s\geq 0$.  It is further clear that $\{\Lambda_{t}\}_{t\geq 0}$ 
forms a semigroup of operators on these spaces. Lemma \ref{wm7} says that
$\calD=C_c^\infty([0,\infty);E)$ is dense in all the spaces above. 
The left translations act strongly continuous on $\calD$ with the sup-norm, and they thus
 form $C_0$-semigroups on $W_{p,\mu}^s(\R_+;E)$ and $H_{p,\mu}^s(\R_+;E)$.

\textbf{(II)} We denote the generator of $\{\Lambda_t\}$ on $W_{p,\mu}^s(\R_+;E)$ by $A$.
We first let $s=k\in \N_0$. To show $A\subseteq \partial_t$, we take $u\in D(A)\subset 
W_{p,\mu}^k(\R_+;E)$. 
Then $u^{(j)}\in L_{1,\text{loc}}([0,\infty);E)$ for all $j\in\{0,\cdots, k\}$
 by Lemma \ref{wm14}, and for  $a,b\in \R_+$ with  $a<b$ it holds 
$$\int_a^b \frac{1}{h}\big (u^{(j)}(\tau+h)-u^{(k)}(\tau)\big )\D \tau 
= \frac{1}{h} \int_b^{b+h} u^{(j)}(\tau) \D \tau -  \frac{1}{h} \int_a^{a+h} u^{(j)}(\tau)\D\tau.
$$ 
As $h\ra 0$, the right hand side converges to $u^{(j)}(b)-u^{(j)}(a)$ for almost all 
$a,b\in \R_+$. The integrand on the left hand side tends  to $A u^{(j)}$ in 
$L_{p,\mu}(\R_+;E)$, and thus in $L_1(a,b; E)$. 
Hence, the left hand side converges to $\int_a^b A u^{(j)}(\tau) \D \tau$. We infer that 
$u\in W_{1,\text{loc}}^{k+1}(\R_+;E)$, with $u^{(j+1)} = Au^{(j)}$, so that
$D(A) \subset W_{p,\mu}^{k+1}(\R_+;E)$ and $\partial_t|_{D(A)} = A$.  

Since $A$ generates a contraction semigroup, we know that
 $1-A$ is invertible. It is further easy to see that $1-\partial_t$ is injective on 
$W_{p,\mu}^{k+1}(\R_+;E)$.  Consequently,  $1- \partial_t= 1-A$, which yields $\partial_t=A$.
This fact implies that $1-\partial_t$ is an isomorphism from $W_{p,\mu}^{k+1}(\R_+;E)$
to $W_{p,\mu}^{k}(\R_+;E)$. By interpolation, we can replace here $k\in\N_0$ 
by any $s\ge0$, and $W$ by $H$.

Finally, let $s\in (k,k+1)$. The first part and Proposition II.2.3
of \cite{EN00} imply that $D(A)= \{u\in W_{p,\mu}^{k+1}(\R_+;E)\,:\, 
u'\in W_{p,\mu}^{s}(\R_+;E)\}$ and $Au=u'$. Using that $1-\partial_t$ is an isomorphism
from $W_{p,\mu}^{s+1}(\R_+;E)$ to $W_{p,\mu}^{s}(\R_+;E)$, 
we conclude that $D(A)=W_{p,\mu}^{s+1}(\R_+;E)$. The same arguments work with $W$ replaced
 by $H$.
 \eprf

Since $W_{p,\mu}^1(\R_+;E)$ is the domain of the generator of  $\{\Lambda_{t}^{E}\}_{t\geq 0}$
by the above lemma, standard interpolation theory yields
$$|u|_{W_{p,\mu}^s( \R_+ ;E)}^p \sim |u|_{L_{p,\mu}(\R_+;E)}^p 
+ \int_0^\infty  t^{-sp-1}\,| \Lambda_t^Eu-u|_{L_{p,\mu}(\R_+;E)}^p \D t$$
for $s\in (0,1)$,  see Proposition~5.7 of \cite{Lun09}. 
Substitution and localization then imply that
\begin{align}
|u|_{W_{p,\mu}^s(J;E)} &\sim |u|_{L_{p,\mu}(J;E)} + [u]_{W_{p,\mu}^s(J;E)} 
 \qquad\text{with} \notag \\
\label{pj956}
[u]_{W_{p,\mu}^s(J;E)}^p &:= \int_0^T\int_0^t \tau^{p(1-\mu)} 
          \frac{|u(t)-u(\tau)|_E^p}{(t-\tau)^{1+sp}} \D \tau\D t
\end{align}
for finite or infinite $J=(0,T)$, $s\in (0,1)$ and $u\in W_{p,\mu}^s(J ;E)$, see Proposition~1.1.13 of \cite{Mey10}.


Using the generation property of $\partial_t$ and the transference principle, 
we show that  $-\partial_t$ with domain  
$W_{p,\mu}^1( \R_+ ;E)$ possesses a bounded $\calH^\infty$-calculus on $L_{p,\mu}(\R_+;E)$. 
Here it is for the first time essential that $E$ is of class $\calH\calT$. As shown in 
\cite{PS04}, the realization of $\partial_t$ with domain ${}_0W_{p,\mu}^1(\R_+;E)$ 
also admits a bounded $\calH^\infty$-calculus on $L_{p,\mu}(\R_+;E)$ (and it is in particular
sectorial), although 
$-\partial_t$  with domain ${}_0W_{p,\mu}^1(\R_+;E)$ does not generate a semigroup 
on $L_{p,\mu}(\R_+;E)$.

\begin{thm}\label{wm11} \textsl{Let $0< 1/p<\mu\leq 1$. Then the 
operators 
$\partial_t$ with domain ${}_0W_{p,\mu}^1(\R_+;E)$ and  $-\partial_t$ with domain
 $W_{p,\mu}^1(\R_+;E)$ possess a bounded $\calH^\infty$-calculus on $L_{p,\mu}(\R_+;E)$
with $\calH^\infty$-angle $\pi/2$. }
\end{thm}
\bprf 
The assertion on $\partial_t$ is proved in Theorem 4.5 of \cite{PS04}. 
To treat the operator $-\partial_t$, we employ the vector-valued transference principle
from \cite{HP98}. We first introduce vector-valued extensions of operators. 
Let $(\Omega,\nu)$ be a measure space and $S$ be a bounded, positive operator on 
$L_p(\Omega, \nu)$. For simple functions $u$ of the form $u = \sum_{i=1}^N u_i x_i$
with $x_i \in E$ and simple $u_i:\Omega\to \C$, we define $S^Eu(\cdot) 
:= \sum_{i=1}^N (Su_i)(\cdot) x_i.$ The operator $S^E$ extends uniquely to $L_p(\Omega, \nu; E)$ 
 with $\|S^E\|_{\calB(L_p(\Omega, \nu; E))} = \|S\|_{\calB(L_p(\Omega, \nu))}$, 
see  Lemma~10.14 of \cite{KW04}.

We consider  the left translation $\Lambda_t^E$, $t\ge0$, on $L_{p,\mu}(\R_+;E) 
= L_p(\R_+, t^{p(1-\mu)}\D t;E)$. Clearly, $\Lambda_t^{E} = (\Lambda_t^{\C})^E$ is the 
vector-valued extension of the scalar left translation $\Lambda_t^{\C}$.
 Due to Lemma~\ref{sec:timeder}, the family $\{ \Lambda_t^{\C}\}_{t\geq 0}$  forms a 
$C_0$-semigroup of positive contractions on $L_{p,\mu}(\R_+)$, and $\partial_t$ with domain
 $W_{p,\mu}^1(\R_+;E)$ generates the  vector-valued extension $\{ \Lambda_t^{E}\}_{t\geq 0}$
on $L_{p,\mu}(\R_+;E)$. Moreover, $\partial_t $ is injective on this space. Theorem~6 of 
\cite{HP98} now yields that $-\partial_t$ admits a bounded $\calH^\infty$-calculus with 
angle $\pi/2$.\eprf

Theorem~\ref{wm11}  implies that  $1+\partial_t: {}_0W_{p,\mu}^{1}(\R_+;E) \to L_{p,\mu}(\R_+;E)$
is isomorphic. The restriction of this  operator to finite $J$ is still injective,
 and its bounded inverse is given by the restriction of the inverse for $J=\R_+$.
By induction and interpolation we deduce that 
 $1+\partial_t: {}_0W_{p,\mu}^{s+1}(J;E) \to {}_0W_{p,\mu}^s(J;E)$
is isomorphic for all $s\ge0$ and $J$. For  $s= [s] + s_*$ with $[s]\in \N_0$
and $s_*\in [0,1)$, we thus obtain the important (and expected) equality 
\begin{align}
\label{lg56}
{}_0W_{p,\mu}^s(J;E) &= \big\{u\in {}_0W_{p,\mu}^{[s]}(J;E)\,:\, 
    u^{([s])}\in {}_0W_{p,\mu}^{s_*}(J;E)\big \},
\end{align}
where the natural norms are equivalent with constants  independent of $J$.
In Lemma~\ref{sec:timeder} we have seen  that $1-\partial_t$ is an isomorphism
from $W_{p,\mu}^{s+1}(\R_+;E)$ to $W_{p,\mu}^s(\R_+;E)$ for all $s\ge0$.
For $J=\R_+$, it follows that
\begin{align}\label{lg55}
W_{p,\mu}^s(J;E) &= \big \{u\in W_{p,\mu}^{[s]}(J;E)\,:\, u^{([s])}\in W_{p,\mu}^{s_*}(J;E)\big\}.
\end{align}
By means of extension and restriction, one can extend this identity to finite $J$.
(Observe that the formula for $\calE_J$ given in the proof of 
Lemma~\ref{sec:ext3} imply that if $u$ is contained the space on right hand side,
 then its extension belongs to the analogous space with $J=\R_+$.)
These characterizations remain valid if one replaces the $W$-spaces by the $H$-spaces.
We next prove general interpolation properties of the weighted spaces.

\begin{lem}\label{lg54}
\textsl{Let $J=(0,T)$ be finite or infinite, $p\in (1,\infty)$, $\mu\in (1/p,1]$, 
$0\leq s_1 < s_2$ and $\theta\in (0,1)$. Set $s= (1-\theta)s_1 + \theta s_2$. 
Then it holds
\beq\label{inter1}
\big [H_{p,\mu}^{s_1}(J;E), H_{p,\mu}^{s_2}(J;E)\big]_\theta = H_{p,\mu}^s(J;E).
\eeq
If also $s\notin \N$, we have
\beq\label{lg59}
\big(H_{p,\mu}^{s_1}(J;E), H_{p,\mu}^{s_2}(J;E)\big)_{\theta,p} = W_{p,\mu}^s(J;E).
\eeq
Moreover, if $s_1, s_2,s\notin \N_0$ then
\begin{align}\label{inter2}
\big[W_{p,\mu}^{s_1}(J;E), W_{p,\mu}^{s_2}(J;E)\big]_{\theta} = W_{p,\mu}^s(J;E), \quad
\big(W_{p,\mu}^{s_1}(J;E), W_{p,\mu}^{s_2}(J;E)\big)_{\theta,p} = W_{p,\mu}^s(J;E). 
\end{align}
If $F\stackrel{d}{\hra} E$ is a Banach space of class $\HT$, then it holds for $\tau\geq 0$ that 
\begin{align}
\big(H_{p,\mu}^\tau(J; E),H_{p,\mu}^\tau(J;F)\big)_{\theta,p} 
  &= H_{p,\mu}^\tau\big(J;(E,F)_{\theta,p}\big),\label{inter3}\\
\big[H_{p,\mu}^\tau(J; E),H_{p,\mu}^\tau(J;F)\big]_{\theta} 
  &=  H_{p,\mu}^\tau\big(J;[E,F]_{\theta}\big). \label{inter4}
  \end{align}
All these assertions remain true if one replaces the $W$- and $H$-spaces by
${}_0W$- and ${}_0H$-spaces, respectively, where the constants
of the norm equivalences do not depend on $T$ if  $s_2\leq 2$.}
\end{lem}
\bprf 
\textbf{(I)} First, let $J=\R_+$. Throughout,
we consider the operator $A = 1-\partial_t$ on $L_{p,\mu}(\R_+;Z)$
with $D(A)=W^1_{p,\mu}(\R_+;Z)$ for any Banach space $Z$ of class $\calH\calT$, 
and we put $ X:= L_{p,\mu}(\R_+;E).$  By Theorem \ref{wm11}, the operator $A$ possesses
a bounded $\calH^\infty$-calculus on $X$ so that $A$ admits bounded 
imaginary powers (see (2.15) of \cite{DHP03}). A theorem by Yagi now implies that
$D(A^\alpha) = [X, D(A)]_\alpha=H_{p,\mu}^\alpha(\R_+;E)$ for all $\alpha\in (0,1)$, see
Theorem~1.15.3 of \cite{Tri94}. For $\alpha\ge1$, formula (\ref{lg55}) further yields
\begin{align}
D(A^\alpha) &\,= \{u\in D(A^{[\alpha]}) \;:\; A^{[\alpha]} u \in D(A^{\alpha-[\alpha]})\} \notag\\
&\, = \{u\in H_{p,\mu}^{[\alpha]}(\R_+;E) \;:\; u^{([\alpha])} \in 
      H_{p,\mu}^{\alpha- [\alpha]}(\R_+;E)\}= H_{p,\mu}^\alpha(\R_+;E). \label{Aalpha}
\end{align}
Employing this equation and again Yagi's theorem, we obtain \eqref{inter1} by computing
\begin{align*}
\big[H_{p,\mu}^{s_1}(\R_+;E), H_{p,\mu}^{s_2}(\R_+;E)\big]_\theta 
  = \big[D(A^{s_1}), D(A^{s_2})\big]_{\theta}= D(A^s) = H_{p,\mu}^s(\R_+;E).
\end{align*}

\textbf{(II)} We continue with (\ref{lg59}). The operator $A^{s_1}$ induces an isomorphism 
$$\big(H_{p,\mu}^{s_1}(J;E), H_{p,\mu}^{s_2}(J;E)\big)_{\theta,p} 
  = \big(D(A^{s_1}), D(A^{s_2})\big)_{\theta,p} \ra \big(X, D(A^\tau)\big)_{\theta,p},$$ 
where $\tau = s_2-s_1$. We first assume that $\theta\tau\notin\N_0$, and put $k=[\theta\tau]\in \N_0$.
It follows from reiteration (see e.g.\ Theorem~1.10.3.2 of \cite{Tri94}) that 
$$\big(X, D(A^\tau)\big)_{\theta,p} = \big(D(A^{k}), D(A^\tau)\big)_{\sigma,p},$$ 
with $\sigma = \frac{\tau\theta - k}{\tau - k}\in (0,1),$ and that the operator 
$A^{k}$ induces an isomorphism 
$$\big(D(A^{k}), D(A^{\tau})\big)_{\sigma,p} \ra \big(X, D(A^{\tau- k})\big)_{\sigma,p}
  = \big(X, D(A)\big)_{\sigma(\tau-k),p}= W_{p,\mu}^{\sigma(\tau-k)}(\R_+;E),$$ 
where $\sigma(\tau-k) = \tau\theta -k \in (0,1)$. On the other hand,
 $A^{-(s_1+k)}$ induces an isomorphism 
$$W_{p,\mu}^{\tau\theta -k}(\R_+;E) 
         \ra W_{p,\mu}^{(s_2-s_1)\theta + s_1}(\R_+;E)=W_{p,\mu}^{s}(\R_+;E).$$ 
Hence, (\ref{lg59}) holds if $\theta(s_2-s_1)\notin\N_0$. 

\textbf{(III)} Next we derive \eqref{inter2}. We take an integer $k>s_2$. We can write 
$W_{p,\mu}^{s_j}(\R_+;E)=\big ( X, H_{p,\mu}^k(\R_+;E)\big )_{s_j/k,p}$ thanks to part (II) and 
$ k \cdot s_j/k =s_j\notin \N_0$. Now reiteration yields 
\[ \big( W_{p,\mu}^{s_1}(J;E), W_{p,\mu}^{s_2}(J;E)\big)_{\theta,p}
  = \big( X, H_{p,\mu}^k(\R_+;E) \big)_{s/k,p} =W_{p,\mu}^s(\R_+;E), \]
  using again part (II) and that $k\cdot s/k \notin\N_0$. The other equation in
  \eqref{inter2} is shown similarly, employing Remark~1.10.3.2 of \cite{Tri94}.
The remaining case $\theta(s_2-s_1)\in\N_0$ in Step II can now be treated by an reiteration 
argument with exponents $\theta\pm \eps$.

\textbf{(IV)} Now let $F\stackrel{d}{\hra} E$ be a Banach space of class $\HT$ 
and $\tau \geq 0$. Due to \eqref{Aalpha}, the operator $A^\tau$ is an isomorphism 
$$\big (H_{p,\mu}^\tau(\R_+; E),H_{p,\mu}^\tau(\R_+;F)\big)_{\theta,p} 
 \ra \big(L_{p,\mu}(\R_+; E),L_{p,\mu}(\R_+;F)\big)_{\theta,p}.$$ 
Theorem 1.18.4 of \cite{Tri94} says that the latter space equals 
$L_{p,\mu}\big(\R_+; (E,F)_{\theta,p}\big)$. Since $A^{-\tau}$ maps this space isomorphically to 
$H_{p,\mu}^\tau\big(\R_+; (E,F)_{\theta,p}\big)$, we have shown \eqref{inter3}.
The formula \eqref{inter4} is shown in the same way.

\textbf{(V)} Replacing the operator $A = 1-\partial_t$ by $A_0:= 1+\partial_t$, the same 
arguments 
as above show the asserted equalities for the ${}_0W$- and the ${}_0H$-spaces. This finishes the 
case $J=\R_+$. The case of a finite interval is then deduced from the half-line case, using the 
extension operators $\calE_J$ and $\calE_J^0$ from Lemma \ref{sec:ext3} and Theorem~1.2.4
of \cite{Tri94} (see also Section~I.2.3 of \cite{Ama95}). The dependence of the norm 
equivalence constants on the length of $J$ carries over from the properties of the extension 
operators. \eprf

We now extend Theorem~\ref{wm11} to fractional powers of the  derivative 
acting on the $W$- and $H$-scales. To this aim,  we  note that if a sectorial 
operator $A$ has a bounded $\calH^\infty$-calculus with $\calH^\infty$-angle
$\phi\in [0,\pi)$ and if $\alpha\in(0,\pi/\phi)$, then $A^\alpha$ has the same property 
with angle less or equal $\alpha\phi$. This fact easily follows from the properties
of the calculus, see Lemma~A.3.5 in \cite{Mey10}.

\begin{prop}\label{maal1}
\textsl{Let $0<1/p<\mu\le 1$, $s\geq 0$, $\alpha\in (0,2)$ and $\omega>0$. 
Then the operators
\begin{align*}
 (\omega-\partial_t)^\alpha &\;\text{ on } \;H_{p,\mu}^s(\R_+;E) 
           \qquad \text{ with domain }\; H_{p,\mu}^{s+\alpha}(\R_+;E),\\
(\omega-\partial_t)^\alpha &\;\text{ on } \;W_{p,\mu}^s(\R_+;E)  
     \qquad \text{ with domain }\; W_{p,\mu}^{s+\alpha}(\R_+;E),\qquad s,s+\alpha\notin \N_0,\\
(\omega+\partial_t)^\alpha &\;\text{ on } \;{}_0 H_{p,\mu}^s(\R_+;E) 
   \qquad \text{ with domain }\; {}_0H_{p,\mu}^{s+\alpha}(\R_+;E),\\
(\omega+\partial_t)^\alpha &\;\text{ on } \;{}_0W_{p,\mu}^s(\R_+;E) 
  \qquad \text{ with domain }\; {}_0W_{p,\mu}^{s+\alpha}(\R_+;E),\qquad s,s+\alpha\notin \N_0,
\end{align*}
are invertible and possess a bounded $\calH^\infty$-calculus with $\calH^\infty$-angle 
less or equal $\alpha\pi/2$.}
\end{prop}
\bprf 
We first consider the case $s=0$. Proposition~2.11 of \cite{DHP03}
and Theorem \ref{wm11} imply that the realization of  $\omega-\partial_t$ on 
$L_{p,\mu}(\R_+;E)$ with domain $W_{p,\mu}^1(\R_+;E)$ admits a bounded $\calH^\infty$-calculus 
with $\calH^\infty$-angle equal to $\pi/2$. As noted above, also  $(\omega-\partial_t)^\alpha$ 
possesses a bounded $\calH^\infty$-calculus with $\calH^\infty$-angle less or equal $\alpha\pi/2$, 
because of $\alpha\in (0,2)$. As in \eqref{Aalpha} we see that then
$D((\omega-\partial_t)^\alpha) =  H_{p,\mu}^\alpha(\R_+;E)$ holds for all $\alpha\ge0$.

Since $\omega-\partial_t$ is invertible and 
$H_{p,\mu}^{\tau+s}(\R_+;E)=D((\omega-A)^{\tau+s})$ holds for all $s,\tau\ge0$, the operator 
$(\omega-\partial_t)^s$ induces an isomorphism from $H_{p,\mu}^{\tau+s}(\R_+;E)$ to 
$H_{p,\mu}^{\tau}(\R_+;E)$. Using that $(\omega-\partial_t)^{-s}$  and 
$(\omega-\partial_t)^\alpha$ commute on $H_{p,\mu}^\alpha(\R_+;E)$, we derive
from Proposition~2.11 of \cite{DHP03} that  $(\omega-\partial_t)^\alpha$ has a 
bounded $\calH^\infty$-calculus on $H_{p,\mu}^{s}(\R_+;E)$ with angle not larger than 
$\alpha\pi/2$, and that its domain equals $H_{p,\mu}^{s+\alpha}(\R_+;E)$.

Let $s, s+\alpha\notin \N_0$. Interpolation of the $H$-case and 
Lemma~\ref{lg54} show that $(\omega-\partial_t)^\alpha$  with domain 
$W_{p,\mu}^{s+\alpha}(\R_+;E)$ has a bounded $\calH^\infty$-calculus 
on $W_{p,\mu}^{s}(\R_+;E)$ with $\calH^\infty$-angle less or equal $\alpha\pi/2$. 
The operator $\omega+\partial_t$ can be treated in the same way.\eprf

The ${}_0W_{p,\mu}^s$-spaces can be characterized in terms of the kernel of the temporal trace, 
as the next proposition shows. These results are mainly due to Grisvard \cite{Gri63}. 
Observe that the limit number for the existence of a trace is $s=1-\mu +1/p$. Therefore, 
if $\mu$ runs through the interval $(1/p,1]$ this limit number runs through the interval 
$[1/p,1)$. Of course, for $\mu=1$ the limit number $s=1/p$ for the unweighted case is recovered. 
We do not treat the limit cases $s= k + 1-\mu+1/p$ with $k\in \N_0$, see Remarque~4.2
of \cite{Gri63} and Remark~3.6.3.2 of \cite{Tri94} for a discussion. We also do not consider the corresponding 
characterizations of the ${}_0H_{p,\mu}^s$-spaces, which are not needed below. They should be 
correct, but it seems that their proofs require a much greater effort.


\begin{prop}\label{wm10}\textsl{Let $J=(0,T)$ be finite or infinite,
$p\in (1,\infty)$ and $\mu \in (1/p,1]$. Then for $0 < s < 1-\mu+1/p$ it holds 
\beq\label{gris1}
C_c^\infty(\ol{J}\bs\{0\};E) \stackrel{d}{\hra} W_{p,\mu}^s(J;E) \quad \text{ and } \quad 
W_{p,\mu}^s(J;E) = {}_0W_{p,\mu}^s(J;E).
\eeq
For $k + 1-\mu+1/p < s < k+1 + (1-\mu+1/p)$ with $k\in \N_0$ we have
\beq\label{kv5001}
W_{p,\mu}^s(J;E) \hra BU\!C^k(\ol{J};E),
\eeq
where here one may replace $W_{p,\mu}^s$ by $H_{p,\mu}^s$, and moreover 
\beq\label{lg62}
{}_0W_{p,\mu}^s(J;E) = \big \{u\in W_{p,\mu}^s(J;E)\;:\; u^{(j)}(0) = 0,\;\; j\in \{0,...,k\}\big \}.
\eeq 
If in addition  $s\in [0,2]$, then  the embedding constants for  
$${}_0W_{p,\mu}^s(J;E) \hra BU\!C^k(\ol{J};E) \quad 
\text{ and } \quad {}_0H_{p,\mu}^s(J;E) \hra BU\!C^k(\ol{J};E)$$ 
do not depend on $J$.}
\end{prop}
\bprf 
In \cite{Gri63} much of the proposition is proved for the $W$-spaces in  
the scalar-valued case $E=\C$ with $J=\R_+$. 
An inspection of the proofs given there shows that they only make use of basic facts, 
interpolation theory and 
the scalar versions of  Lemmas~\ref{wm7} and \ref{sec:timeder} above. 
Thus the results of \cite{Gri63} carry over to the case of a general Banach space $E$. 
Moreover, the case of a finite interval is obtained from the half-line  case by extension 
and restriction. One can replace $W$ by $H$ as asserted because of
(\ref{kv33}) and (\ref{kv5000}).

Assertion \eqref{gris1} is shown in Th\'eor\`eme~2.1 and  Th\'eor\`eme~4.1 of \cite{Gri63}. The 
embedding \eqref{kv5001} is  proved for $k=0$ in Th\'eor\`eme~5.2 of \cite{Gri63}, 
and the general case $k\in \N$ is an immediate consequence in view of \eqref{lg55}. 
For $s\leq 1$,  identity (\ref{lg62}) is shown in Proposition~1.2 and Th\'eor\`eme~4.1 of
\cite{Gri63}. For integers $s=k+1$ this equality holds by definition. The general case then 
follows from (\ref{lg56}),  \eqref{gris1} and (\ref{lg62}) for $s-[s]$, as well as (\ref{lg55}). 
\eprf

We next establish embeddings of Sobolev type for $W_{p,\mu}^s$ and $H_{p,\mu}^s$.

\begin{prop}\label{sec:rol2}\textsl{Let $J=(0,T)$ be finite, $1<p<q<\infty$, $\mu>1/p$
and $s>\tau \geq 0$. Then 
\begin{align}\label{lg63}
W_{p,\mu}^s(J;E) &\hra W_{q,\mu}^\tau(J;E) \qquad  \text{holds if \ }\;s-(1-\mu+1/p) 
       > \tau - \frac{p(1-\mu+1/p)}{q},\\
\label{lg64}
W_{p,\mu}^s(J;E) &\hra W_{q}^\tau(J;E) \qquad \text{holds if \ }\;s - (1-\mu+1/p) > \tau - 1/q.
\end{align}
These embeddings remain true if one replaces the $W$-spaces by the $H$-, the ${}_0W$- and the 
${}_0H$-spaces, respectively. In the two latter cases, restricting to $s\in [0,2]$, 
for given $T_0>0$ the embeddings hold with a uniform constant for all $0< T\leq T_0$.}
\end{prop}

\bprf Throughout this proof, let $T_0>0$ be given. Since the inequality signs in (\ref{lg63}) 
and (\ref{lg64}) are strict, we may assume that $s\notin \N$. Again we only have to consider 
the $W$-case due to (\ref{kv33}) and (\ref{kv5000}).

\textbf{(I)} We prove (\ref{lg63}) for $\tau=0$. For $s>1-\mu+1/p$,  Proposition \ref{wm10} 
yields 
$$W_{p,\mu}^s(J;E) \hra L_\infty(J;E)\hra L_{q,\mu}(J;E) \qquad \text{for all }\;q\in (p,\infty),$$ 
with the asserted behaviour of the embedding constant in the ${}_0W$-case. 
If $s\leq 1-\mu+1/p$, take $q\in (p,\infty)$
as in  \eqref{lg63}. Choose  $\eta\in (1-\mu+1/p,1)$ and  $r\in (p,\infty)$
such that  $\frac{1}{q} \le \frac{1-s/\eta}{p} + \frac{s/\eta}{r}$. We have already shown that
$W_{p,\mu}^\eta(J;E) \hra L_{r,\mu}(J;E)$.  Lemma~\ref{lg54} and 
Theorem~1.18.5 of \cite{Tri94} then imply 
$$W_{p,\mu}^s(J;E) = \big (L_{p,\mu}(J;E),W_{p,\mu}^\eta(J;E)\big)_{s/\eta,p} 
     \hra \big(L_{p,\mu}(J;E),L_{r,\mu}(J;E)\big)_{s/\eta,p} \hra L_{q,\mu}(J;E).$$ 
In the ${}_0W$-case, all embedding constants are uniform in $T\leq T_0$ for $s\in [0,2]$.

\textbf{(II)} To prove (\ref{lg63}) for $\tau>0$, we set 
$\alpha= (1-p/q) (1-\mu+1/p)>0$ for an exponent $q>p$ as in (\ref{lg63}). We fix some
$\eps>0$ and $k\in\N_0$ such that $s-\eps>\tau+\alpha>0$,  $k+\alpha< s-\eps<k+1+\alpha$, 
and $\kappa:= k+\alpha+\eps$ is not an integer.
 Lemma \ref{lg54} now yields
$$W_{p,\mu}^s(J;E) = \big(W_{p,\mu}^\kappa(J;E), W_{p,\mu}^{\kappa+1}(J;E)\big)_{s-\kappa, p},$$
 Using (\ref{lg55}) and (\ref{lg63}) with $\tau = 0$, we obtain  
$$W_{p,\mu}^\kappa(J;E) \hra W_{q,\mu}^k(J;E) \quad \text{ and } \quad
  W_{p,\mu}^{\kappa+1}(J;E) \hra W_{q,\mu}^{k+1}(J;E).$$
Since $(\cdot, \cdot)_{\theta,p} \hra (\cdot, \cdot)_{\theta,q}$ for $\theta\in (0,1)$, it 
follows
$$W_{p,\mu}^s(J;E) \hra \big(W_{q,\mu}^k(J;E), W_{q,\mu}^{k+1}(J;E)\big)_{s-\kappa, q} 
= W_{q,\mu}^{k + s-\kappa}(J;E) \hra W_{q,\mu}^{\tau}(J;E).$$

\textbf{(III)} For ${}_0W$-spaces, the dependence on $T$ for $s\in [0,2]$ carries over from 
Lemma~\ref{lg54} and (\ref{lg63}) with $\tau = 0$. 
The embedding (\ref{lg64}) is shown similarly employing Lemma~\ref{wm14},
see Proposition~1.1.12 in \cite{Mey10}.\eprf

We finish this section with Poincar\'e's inequality in the weighted spaces. 

\begin{lem}\label{sec:poincare}\textsl{Let $J=(0,T)$ be finite,
$p\in (1,\infty)$, $\mu \in (1/p,1]$ and $s\in [0,1)$. Then it holds
$$
|u|_{L_{p,\mu}(J;E)} \lesssim T\,|u'|_{L_{p,\mu}(J;E)}, \qquad |u|_{{}_0W_{p,\mu}^s(J;E)} + |u|_{{}_0H_{p,\mu}^s(J;E)} \lesssim T^{1-s}\,|u|_{W_{p,\mu}^1(J;E)},
$$
for all $u\in {}_0W_{p,\mu}^1(J;E)$.}
\end{lem}
\bprf 
Let $u\in {}_0W_{p,\mu}^1(J;E).$ Using H\"older's inequality, we estimate
$$
t^{p(1-\mu)}|u(t)|_E^p  \leq   
         t^{p(1-\mu)} \left ( \int_{0}^{T} s^{-(1-\mu)} s^{1-\mu}|u'(s)|_E \D s \right )^p
  \lesssim t^{p(1-\mu)} \, T^{(1-p'(1-\mu))p/p'}  |u'|_{L_{p,\mu}(J;E)}^p
$$
for $t\in J$. The first assertion now follows by integration over $J$. For $s\in [0,1)$ we have
$$|u|_{{}_0W_{p,\mu}^s(J;E)} + |u|_{{}_0H_{p,\mu}^s(J;E)} \lesssim  
    |u|_{{}_0W_{p,\mu}^1(J;E)}^s|u|_{L_{p,\mu}(J;E)}^{1-s}$$ 
by interpolation, implying the second asserted inequality. \eprf

\section{Weighted anisotropic spaces}\label{was}
Let again $E$ be a Banach space of class $\HT$, let $J=(0,T)$ be finite or infinite, and 
let further $\Omega\subset \R^n$ be a domain with compact smooth boundary $\partial\Omega$, 
or $\Omega\in \{\R^n, \R_+^n\}$. For $p,q\in (1,\infty)$ and $r>0$ we denote by 
$$H_p^r(\Omega;E), \qquad W_p^r(\Omega;E), \qquad B_{p,q}^r(\Omega;E),$$ 
the $E$-valued Sobolev, Slobodetskii and Besov spaces, respectively. 
If $\Omega=\R^n$ we also use this spaces with $r\le 0$. For the definitions and 
properties of these spaces we refer to \cite{Ama97},  \cite{Sch} or \cite{Zi89}. 
The scalar-valued case is treated extensively in \cite{Tri94}. In particular, we have 
$$B_{p,p}^r(\Omega;E)= W_p^r(\Omega;E) \qquad \text{for all  \ } p\in [1,\infty), \quad 
      r\notin \N_0.$$ 
The corresponding spaces over the boundary $\partial\Omega$ of $\Omega$ are defined via local 
charts as in Definition~3.6.1 of \cite{Tri94}, for instance.

We mainly investigate weighted anisotropic spaces, i.e., intersections of spaces of the form
\beq\label{aniso}
H_{p,\mu}^{s}\big (J; H_p^r(\Omega;E)\big), \ \ \  W_{p,\mu}^{s}\big(J; W_p^r(\Omega;E)\big), 
\ \ \  H_{p,\mu}^{s}\big(J; W_p^r(\Omega;E)\big), \ \ \ 
W_{p,\mu}^{s}\big(J; W_p^r(\Omega;E)\big),
\eeq
where $s,r\geq 0$. In what follows we refer to $t\in J$ as time and to $x\in \ol{\Omega}$ as 
space variables.

We start with two important tools. First, given $k\in \N$, there is an extension operator 
$\calE_\Omega$ to $\R^n$ for functions defined on $\Omega$ (i.e., we have
$(\calE_\Omega u)|_{\Omega} = u$) satisfying
\beq\label{ext1}
\calE_\Omega\in \calB\big ( B_{p,q}^r(\Omega;E), B_{p,q}^r(\R^n;E)\big )\cap 
 \calB\big ( H_p^r(\Omega;E), H_p^r(\R^n;E)\big )
\eeq
 for all $p,q\in (1,\infty)$ and $r\in [0,k]$.
For integer $r\in [0,k]$, the constructions in the proofs of Theorems~5.21 and  5.22 \cite{AF03} 
for the scalar-valued spaces literally carries over to the vector-valued case. 
The general case $r\in [0,k]$ can be treated via interpolation. 
Applying $\calE_\Omega$ pointwise almost everywhere in time and using again interpolation, 
we obtain a spatial extension
 operator for the anisotropic spaces, again denoted by $\calE_\Omega$:
\beq\label{ext2}
\calE_{\Omega}\in \calB\big (H_{p,\mu}^{s}(J; H_p^r(\Omega;E)), H_{p,\mu}^{s}(J; H_p^r(\R^n;E))\big), \qquad s\geq 0, \quad r\in [0,k].
\eeq 
Of course, here a $H$-space may be replaced by a $W$-space at the first or the second or at 
both positions, and this remains true for the ${}_0H_{p,\mu}^s$- and the 
${}_0W_{p,\mu}^s$-spaces with respect to time.

Second, we consider operators build from the time derivative and the Laplacian
whose domains are given by intersections of the spaces in \eqref{aniso} with $\Omega=\R^n$ and $J=\R_+$.
This fact will allow us to study the anisotropic spaces by means of these operators.

\begin{lem}\label{wm22}\textsl{Let $E$ be a Banach space of class $\calH\calT$, let $p\in (1,\infty)$, $\mu\in (1/p,1]$, $s,r\geq 0$, $\alpha\in (0,2)$ and  $\beta >0$, let $\omega,\omega'\geq 0$ satisfy $\omega+\omega'\neq 0$ and set $$H_{p,\mu}^s(H_p^r) := H_{p,\mu}^s\big (\R_+;H_p^r(\R^n;E)\big), \qquad {}_0H_{p,\mu}^s(H_p^r):= {}_0H_{p,\mu}^s\big(\R_+;H_p^r(\R^n;E)\big),$$ and analogously for the other types of spaces in (\ref{aniso}). Let $\Delta_n$ be the Laplacian on $\R^n$. Then the following holds true. 
\begin{itemize}
\item[\textbf{a)}] The pointwise realization of 
$(\omega -\Delta_{n})^{\beta/2}$ on any of the spaces
\begin{align*}
&H_{p,\mu}^s(H_p^r), \quad \text{ with domain }\;\;  H_{p,\mu}^s(H_p^{r+\beta}),\\
 &H_{p,\mu}^s(W_p^r), \quad \text{ with domain }\;\;  H_{p,\mu}^s(W_p^{r+\beta}), \quad r,r+\beta\notin \N_0,\\
 &W_{p,\mu}^s(H_p^r), \quad \text{ with domain } \;\; W_{p,\mu}^s(H_p^{r+\beta}), \\
 &W_{p,\mu}^s(W_p^r), \quad \text{ with domain }\;\;  W_{p,\mu}^s(W_p^{r+\beta}), \quad r,r+\beta\notin \N_0,
\end{align*}
admits a bounded $\calH^\infty$-calculus with $\calH^\infty$-angle equal to zero,
and it is invertible if $\omega>0$. Here one can replace the $H_{p,\mu}^s,W_{p,\mu}^s$-spaces by the ${}_0H_{p,\mu}^s,{}_0W_{p,\mu}^s$-spaces, respectively.
\item[\textbf{b)}] The operator $L:=(\omega' -\partial_t)^{\alpha} + (\omega -\Delta_{n})^{\beta/2}$ acting on any of the spaces
\begin{align*}
 &H_{p,\mu}^s(H_p^r), \quad \text{ with domain }\;\; H_{p,\mu}^{s+\alpha}(H_p^{\beta}) \cap H_{p,\mu}^s(H_p^{r+\beta}),\\
 &H_{p,\mu}^s(W_p^r), \quad \text{ with domain }\;\; H_{p,\mu}^{s+\alpha}(W_p^{\beta}) \cap H_{p,\mu}^s(W_p^{r+\beta}), \quad r,r+\beta\notin \N_0,\\
 &W_{p,\mu}^s(H_p^r), \quad \text{ with domain } \;\;W_{p,\mu}^{s+\alpha}(H_p^{\beta}) \cap W_{p,\mu}^s(H_p^{r+\beta}), \quad s,s+\alpha \notin \N_0,\\
 &W_{p,\mu}^s(W_p^r), \quad \text{ with domain }\;\; W_{p,\mu}^{s+\alpha}(W_p^{\beta}) \cap W_{p,\mu}^s(W_p^{r+\beta}), \quad s,s+\alpha,r,r+\beta\notin \N_0,
\end{align*}
is invertible and admits bounded imaginary powers with power angle not larger than $\alpha \pi/2$. This remains true for the operator $L_0:=(\omega' + \partial_t)^{\alpha} + (\omega -\Delta_{n})^{\beta/2}$ if one replaces the $H_{p,\mu}^s,W_{p,\mu}^s$-spaces by the ${}_0H_{p,\mu}^s,{}_0W_{p,\mu}^s$-spaces, respectively.
\item[\textbf{c)}] For $\tau\in (0,1]$ it holds $$D(L^\tau) = D((\omega' -\partial_t)^{\alpha\tau})\cap D((\omega -\Delta_n)^{\beta\tau/2}),$$ $$D_L(\tau,p) = D_{(\omega' -\partial_t)^{\alpha}}(\tau,p)\cap D_{(\omega -\Delta_{n})^{\beta/2}}(\tau,p),$$ and this remains true if one replaces $L$ by $L_0$ and $\omega' -\partial_t$ by $\omega' +\partial_t$.
\end{itemize}
}
\end{lem}
\bprf Since $E$ is of class $\HT$, the operator $-\Delta_n$ with domain $H_p^2(\R^n;E)$ 
possesses a bounded $\calH^\infty$-calculus  on $L_p(\R^n;E)$
with $\calH^\infty$-angle equal to zero, see  Theorem~5.5 of \cite{DHP03}. 
As in Proposition~\ref{maal1}, the fractional power $(\omega-\Delta_n)^{\beta/2}$  
with domain $H_p^\beta(\R^n;E)$ has the same property. If $\omega>0$, it is invertible
in  $L_p(\R^n;E)$.

Let $r\geq 0$.
Using $(\omega -\Delta_n)^{r/2}$ as an isomorphism between $H_p^r(\R^n;E)$ and $L_p(\R^n;E)$, 
it then follows from Proposition~2.11 of \cite{DHP03} that $(\omega -\Delta_n)^{\beta/2}$ 
 with domain $H_p^{r+\beta}(\R^n;E)$ possesses a bounded $\calH^\infty$-calculus  on 
$H_p^r(\R^n;E)$ with $\calH^\infty$-angle equal to zero. By interpolation, this fact remains 
true if one considers $(\omega-\Delta_n)^{\beta/2}$ with domain 
$W_p^{r+\beta}(\R^n;E)$  in $W_p^r(\R^n;E)$, provided $r,r+\beta \notin \N_0$. 
It is straightforward to see that  
these properties carry over to the pointwise realizations on the anisotropic spaces, cf. Lemma~A.3.6 in \cite{Mey10}.

Since all the spaces under consideration are of class $\HT$, 
Proposition~\ref{maal1} implies that $(\omega'-\partial_t)^\alpha$  with domain 
$H_{p,\mu}^{s+\alpha}(H_p^r)$ admits a bounded $\calH^\infty$-calculus  on $H_{p,\mu}^s(H_p^r)$
 with $\calH^\infty$-angle less or equal $\alpha\pi/2$, and on the corresponding spaces where 
$H$ is replaced by $W$, with the asserted exceptions. Using these facts and that the resolvents 
of $(\omega'-\partial_t)^\alpha$ and $(\omega-\Delta_n)^{\beta/2}$ commute, the assertions 
for $L$ are 
a consequence of Proposition \ref{dv}. The same arguments show the assertions for $L_0$. 
Finally, c) is a consequence of Lemma~9.5 in  \cite{EPS03} and Teorema 5 in \cite{Gri73} (see also Lemma~3.1
in \cite{DHP07}).
\eprf

With the help of the operators from Lemma \ref{wm22} we establish fundamental embeddings for the 
anisotropic spaces with regularity exponents $(s,r)$ for time and space.
Roughly speaking, we prove that the regularities  $(s+\alpha,r)$ and $(s,r+\beta)$ 
imply all regularities on the line segment from  $(s+\alpha,r)$ to $(s,r+\beta)$. 

\begin{prop}\label{embani}\textsl{Let $J=(0,T)$ be finite or infinite, $0<1/p<\mu\le 1$,
and let $\Omega\subset \R^{n}$ be a domain with compact smooth boundary $\partial\Omega$, or $\Omega\in \{\R^n,\R_+^n\}$. Let further $$s,r\geq 0, \qquad \alpha\in (0,2), \qquad \beta>0, \qquad \sigma\in [0,1],$$ and set $H_{p,\mu}^s(H_p^r) := H_{p,\mu}^s\big (J; H_p^r(\Omega;E)\big ),$ and analogously for the other anisotropic spaces. Then it holds
\beq
H_{p,\mu}^{s+\alpha}(H_p^r) \cap H_{p,\mu}^s(H_p^{r+\beta}) \hra H_{p,\mu}^{s+\sigma\alpha}(H_p^{r+ (1-\sigma)\beta}), \label{maal2}
\eeq
and moreover each of the spaces 
$$H_{p,\mu}^{s+\alpha}(W_p^r) \cap H_{p,\mu}^s(W_p^{r+\beta}), \quad W_{p,\mu}^{s+\alpha}(H_p^r) \cap W_{p,\mu}^s(H_p^{r+\beta}), \quad W_{p,\mu}^{s+\alpha}(H_p^r) \cap H_{p,\mu}^s(W_p^{r+\beta}),$$ 
is continuously embedded into 
$$W_{p,\mu}^{s+\sigma\alpha}(H_p^{r+ (1-\sigma)\beta})\cap H_{p,\mu}^{s+\sigma\alpha}(W_p^{r+ (1-\sigma)\beta}),$$ 
provided all the occurring $W_{p,\mu}$- and $W_p$-spaces have a noninteger order of  differentiability.  
Finally, assuming all orders of differentiability to be noninteger, we have
\beq\label{lg69}
W_{p,\mu}^{s+\alpha}(W_p^r) \cap W_{p,\mu}^s(W_p^{r+\beta}) \hra W_{p,\mu}^{s+\sigma\alpha}(W_p^{r+ (1-\sigma)\beta}).
\eeq
These embeddings remain true if one replaces $\Omega$ by its boundary $\partial\Omega$. They are 
also valid if one replaces all  $H_{p,\mu}$-, $W_{p,\mu}$- spaces by ${}_0H_{p,\mu}$-, ${}_0W_{p,\mu}$-spaces, respectively. Restricting in the latter case to $s+\alpha \leq 2$, the embedding constants do not depend on the length of $J$.}
\end{prop}
\bprf Using extensions and restrictions, and employing that the spaces over $\partial\Omega$ are 
defined via local charts, it suffices to consider the case $J\times \Omega=\R_+\times \R^n$. The 
dependence of the embedding constants on $J$ carries over from the properties of the extension 
operators.

\textbf{(I)} For (\ref{maal2}) we consider the operators $(1-\partial_t)^\alpha$ and 
$(1-\Delta_n)^{\beta/2}$ on $H_{p,\mu}^s(H_p^r)$, which were treated in Proposition~\ref{maal1} and 
Lemma~\ref{wm22}. Note that we have to restrict to $\alpha\in (0,2)$ to obtain sectoriality of 
$(1-\partial_t)^\alpha$. Due to the invertibility of these operators
and the description of their domains, the expression 
 $$|(1-\partial_t)^{\alpha\sigma}(1-\Delta_n)^{\beta(1-\sigma)/2} \cdot|_{H_{p,\mu}^s(H_p^r)}$$ 
is an equivalent norm on $H_{p,\mu}^{s+\sigma\alpha}(H_p^{r+ (1-\sigma)\beta})$
for all $\sigma\in (0,1)$. Since the sum $L$ of these operators is invertible 
by Lemma \ref{wm22}, we obtain the equivalent norm
$$|((1-\partial_t)^{\alpha} + (1-\Delta_n)^{\beta/2})\cdot|_{H_{p,\mu}^s(H_p^r)}$$ 
on $D(L)=H_{p,\mu}^{s+\alpha}(H_p^r) \cap H_{p,\mu}^s(H_p^{r+\beta})$. Hence
 (\ref{maal2}) follows directly from Proposition \ref{dv}. The same arguments show 
that  the embeddings
\begin{align*}
H_{p,\mu}^{s+\alpha}(W_p^r) \cap H_{p,\mu}^s(W_p^{r+\beta}) 
&\hra   H_{p,\mu}^{s+\sigma\alpha}(W_p^{r+ (1-\sigma)\beta}),\\
W_{p,\mu}^{s+\alpha}(H_p^r) \cap W_{p,\mu}^s(H_p^{r+\beta}) 
&\hra      W_{p,\mu}^{s+\sigma\alpha}(H_p^{r+ (1-\sigma)\beta}),
\end{align*}
and (\ref{lg69}) hold, with the asserted exceptions.

\textbf{(II)} We derive the remaining embeddings from (\ref{maal2}) by suitable interpolation 
arguments, which were indicated in Remark~5.3 of \cite{EPS03} in a more specific situation. We 
concentrate on the case $W_{p,\mu}^{s+\alpha}(H_p^r) \cap H_{p,\mu}^s(W_p^{r+\beta})$, the other 
assertions are obtained similarly, see Proposition~1.3.12 of \cite{Mey10}.
 For $s+\alpha, r+\beta\notin \N$ and sufficiently small 
$\eps>0$, we apply the real interpolation functor $(\cdot,\cdot)_{1/2,p}$ to the embeddings
\begin{align*}
H_{p,\mu}^{s+\alpha(1\pm\eps/\beta)}(H_p^{r})\cap H_{p,\mu}^s(H_p^{r+\beta\pm\eps}) 
&\hra    H_{p,\mu}^{s+\sigma\alpha}(H_p^{r+ (1-\sigma)\beta\pm \eps}),\\
H_{p,\mu}^{s+\alpha(1\pm\eps/\beta)}(H_p^{r})\cap H_{p,\mu}^s(H_p^{r+\beta\pm\eps}) 
&\hra H_{p,\mu}^{s+\alpha(\sigma\pm\eps/\beta)}(H_p^{r+ (1-\sigma)\beta}).
\end{align*}
By Lemma \ref{lg54} the terms in the first and in the second line on the right-hand side 
interpolate to  $H_{p,\mu}^{s+\sigma\alpha}(W_p^{r+ (1-\sigma)\beta})$ and 
$W_{p,\mu}^{s+\sigma\alpha}(H_p^{r+ (1-\sigma)\beta})$, respectively. To interpolate the 
left-hand side, we consider the operator 
$$L= (1-\partial_t)^{\alpha(1+\eps/\beta)} + (1-\Delta)^{(\beta+\eps)/2},$$ 
on $H_{p,\mu}^s(H_p^r)$
with domain $D(L)=H_{p,\mu}^{s+\alpha(1+\eps/\beta)}(H_p^{r})\cap H_{p,\mu}^s(H_p^{r+\beta+\eps})$. 
Lemma~\ref{wm22} and Proposition~\ref{maal1} imply
$$D(L^{(\beta-\eps)/(\beta+\eps)})= H_{p,\mu}^{s+\alpha(1-\eps/\beta)}(H_p^{r})\cap 
H_{p,\mu}^s(H_p^{r+\beta-\eps}),$$ 
and reiteration yields
$$\big (D(L^{(\beta-\eps)/(\beta+\eps)}), D(L)\big )_{1/2,p} = 
     D_L\big((1+(\beta-\eps)/(\beta+\eps))/2,p\big).$$ 
Using again Lemma~\ref{wm22} and Proposition~\ref{maal1}, we see  that the latter space equals
 $W_{p,\mu}^{s+\alpha}(H_p^r) \cap H_{p,\mu}^s(W_p^{r+\beta})$, as required.

\textbf{(III)} Starting in Step I with $(1+\partial_t)^\alpha$ instead of $(1-\partial_t)^\alpha$, 
the same arguments as above show that the asserted embeddings are also true for the 
${}_0H_{p,\mu}$- and ${}_0W_{p,\mu}$-spaces.
\eprf

\begin{rem}\emph{The proof shows that in the embeddings where only Proposition \ref{dv} was used, the orders of integrability in space and time do not have to coincide. In fact, 
the assertions of Lemma \ref{wm22} remain true for the 
extension of the Laplacian on $H_q^r(\R^n;E)$ to $H_{p,\mu}^s\big(J; H_q^r(\R^n;E)\big)$, where
$p,q\in (1,\infty)$ and $\mu\in (1/p,1]$. As in Step I of the proof above, we then derive
$H_{p,\mu}^{s+\alpha}(H_q^r) \cap H_{p,\mu}^s(H_q^{r+\beta}) 
\hra H_{p,\mu}^{s+\sigma\alpha}(H_q^{r+ (1-\sigma)\beta}).$ \hfill \BlackBox}\end{rem}

A typical application of Proposition~\ref{embani} is the following result on the mapping behavior 
of the spatial derivative on anisotropic $H$-spaces. 

\begin{lem}\label{spatiald} 
\textsl{Let $E$ be a Banach space of class $\HT$, let $J=(0,T)$ be finite or finite, and let 
$\Omega\subset \R^n$ be a domain with compact smooth boundary, or $\Omega\in \{\R_+^n,\R^n\}$. 
Let further 
$$s\geq 0, \qquad r\in [0,1), \qquad \alpha\in (0,2), \qquad \beta\geq 1,
\qquad 0<1/p<\mu\le1.$$ 
Then the pointwise realization of $\partial_{x_i}$, $i\in \{1,...,n\}$,  is a continuous map 
$$ H_{p,\mu}^{s+\alpha}(H_p^r)\cap H_{p,\mu}^{s}(H_p^{r+\beta}) 
\ra H_{p,\mu}^{s+\alpha - \alpha/\beta}(H_p^r)\cap H_{p,\mu}^{s}(H_p^{r+\beta-1}).$$ 
Its operator norm is independent of $T$ if we restrict to $s+\alpha \leq 2$ and 
to ${}_0H_{p,\mu}$-spaces.}
\end{lem}
\bprf 
By extension and restriction it suffices to consider the case $J\times \Omega=\R_+\times \R^n$. 
Clearly the operator $\partial_{x_i}$ maps continuously $H_{p,\mu}^{s+\alpha}(H_p^r) 
\cap H_{p,\mu}^s(H_p^{r+\beta})\ra H_{p,\mu}^s(H_p^{r+\beta-1}).$  Proposition~\ref{embani} 
further yields the embedding 
$$H_{p,\mu}^{s+\alpha}(H_p^r) \cap H_{p,\mu}^s(H_p^{r+\beta}) 
\hra H_{p,\mu}^{s+\alpha - \alpha/\beta}(H_p^{r+1}).$$ 
Hence, $\partial_{x_i}$ is also a  continuous map from $H_{p,\mu}^{s+\alpha}(H_p^r) 
\cap H_{p,\mu}^s(H_p^{r+\beta})$ to $H_{p,\mu}^{s+\alpha - \alpha/\beta}( H_p^r)$.
\eprf

\section{The temporal and the spatial trace}
We first  consider the temporal trace on anisotropic spaces. Using integration by parts, 
one can see that the representation 
\beq\label{lg80}
u(0) = (2-\mu) \Big ( \sigma^{-(2-\mu)} \int_0^\sigma \tau^{1-\mu}u(\tau) \D \tau 
      - (2-\mu) \int_0^\sigma t^{-(3-\mu)} \int_0^t \tau^{1-\mu}(u(t)-u(\tau))\D \tau \D t \Big )
\eeq
holds true for all $\sigma>0$, Banach spaces $X$ and $u\in W_{1,\text{loc}}^1([0,\infty); X)$.
Let $-A$  generate an exponentially stable analytic $C_0$-semigroup on $X$ and $\theta\in (0,1)$. 
Then the norm of the space $D_A(\theta,p)$ is equivalent to the norm given by
\beq\label{lg83}
|x|_{D_A(\theta,p),*}^p = \int_0^\infty \sigma^{p(1-\theta)} |A e^{-\sigma A}x|_X^p 
       \frac{\D \sigma}{\sigma},
\eeq 
see Theorem~1.14.5 of \cite{Tri94}.
The formula (\ref{lg80}) is the key to the following abstract trace theorem. 

\begin{lem}\label{lg72}\textsl{Let $X$ be a Banach space, $p\in (1,\infty)$, $\mu\in (1/p,1]$, and let the operator $A$ on $X$ with domain $D(A)$ be invertible and admit bounded imaginary powers with power angle strictly smaller than $\pi/2$. Let $s\in (0,1-\mu+1/p)$ and $\alpha> 0$ satisfy $s+\alpha \in (1-\mu+1/p,1)$. Then the temporal trace $\emph{\tr}_0$, i.e., $\emph{\tr}_0u = u|_{t=0}$, maps continuously
\beq\label{lg78}
W_{p,\mu}^{s+\alpha}\big (\R_+; D(A^{s})\big)\cap W_{p,\mu}^{s}\big(\R_+;D(A^{s+\alpha})\big) \ra D_A\big(2s + \alpha - (1-\mu+1/p),p\big).
\eeq
Moreover, $\emph{\tr}_0$ is for $\alpha \in ( 1-\mu+1/p, 1]$ continuous
\beq\label{lg76}
W_{p,\mu}^\alpha\big(\R_+; X\big) \cap L_{p,\mu}\big(\R_+; D_A(\alpha,p)\big) \ra D_A\big(\alpha- (1-\mu+1/p),p\big),
\eeq
and for $s\in (0,1-\mu+1/p)$ it is continuous 
\beq\label{lg77}
W_{p,\mu}^1\big(\R_+; D_A(s,p)\big)\cap W_{p,\mu}^s\big(\R_+;D(A)\big) \ra D_A\big(1+s- (1-\mu+1/p),p\big).
\eeq}
\end{lem}
\bprf  
Lemma~11 and 12 in \cite{DiB84} establish the variants of  
the embeddings (\ref{lg76}) and (\ref{lg77}) in the case without weight. The proofs for 
$\mu\in(1/p,1)$ are similar, using (\ref{lg80}), the representation (\ref{pj956}) of the weighted 
Slobodetskii seminorm and Lemma~\ref{sec:Hardy}, and therefore we omit them. Instead, we 
 concentrate on (\ref{lg78}), taking $u$ from the function space on the left-hand side there and assuming that
$s\in (0,1-\mu+1/p)$, $\alpha> 0$ and  $s+\alpha \in (1-\mu+1/p,1)$.
From (\ref{lg80}) and (\ref{lg83})  we deduce 
\begin{align}
|u(0)&|_{D_A(2s + \alpha - (1-\mu+1/p),p),*}^p 
 = \int_0^\infty \sigma^{p(1-(2s+\alpha-(1-\mu+1/p))} |A e^{-\sigma A} u(0)|_X^p \sigma^{-1}\D 
     \sigma \nonumber\\
&\lesssim \int_0^\infty \sigma^{-p(2s+\alpha)} \Big ( \int_0^\sigma \tau^{1-\mu} 
       |A e^{-\sigma A}u(\tau)| \D \tau\Big )^p \D \sigma\label{lg81}\\
& \quad + \int_0^\infty \Big (  \sigma^{2-\mu-(2s+\alpha)} \int_0^\sigma t^{-(3-\mu)} 
  \int_0^t \tau^{1-\mu} |A e^{-\sigma A}(u(t)-u(\tau))| \D \tau \D t\Big )^p \D \sigma.\nonumber
\end{align}
Recall that
\beq\label{lg82}
|Ae^{-\sigma A}x|_X \lesssim \sigma^{-1+\theta} |A^{\theta} x|_X
\eeq
holds for all $\theta\in (0,1)$ and $x\in X$.
Employing H\"older's inequality, (\ref{lg82}), (\ref{pj956}), Lemma~\ref{sec:Hardy} and 
Proposition~\ref{wm10}, we estimate the first summand in (\ref{lg81}) by
\begin{align*}
&\int_0^\infty\,\sigma^{-p(2s+\alpha)} \Big ( \int_0^\sigma \tau^{1-\mu} |A e^{-\sigma A}u(\tau)| \D \tau\Big )^p \D \sigma \\
&\, \leq\int_0^\infty \int_0^\sigma \tau^{p(1-\mu)} |A e^{-\sigma A}u(\tau)|^p \sigma^{p-1} \sigma^{-p(2s+\alpha)} \D \tau \D \sigma\\
&\, \lesssim \int_0^\infty \int_0^\sigma \tau^{p(1-\mu)} |A^{s+\alpha}u(\tau)|^p \sigma^{-(1+ps)} \D \tau \D \sigma\\
&\, \lesssim \int_0^\infty \Big ( \int_0^\sigma \tau^{p(1-\mu)} 
|(u(\sigma)-u(\tau))|_{D(A^{s+\alpha})}^p (\sigma-\tau)^{-(1+ps)} \D \tau 
   +  \sigma^{p(1-\mu -s)} |u(\sigma)|_{D(A^{s+\alpha})}^p \Big ) \D \sigma\\
&\, \lesssim |u|_{W_{p,\mu}^s(\R_+; D(A^{s+\alpha}))}^p.
\end{align*}
We further use (\ref{lg82}),  Hardy's inequality \eqref{hardy-ineq}, H\"older's inequality 
and (\ref{pj956})  to control the second summand in (\ref{lg81}) by
\begin{align*}
 \int_0^\infty  &\,  \Big ( \sigma^{2-\mu-(2s+\alpha)} \int_0^\sigma t^{-(3-\mu)} \int_0^t \tau^{1-\mu} |A e^{-\sigma A}(u(t)-u(\tau))| \D \tau \D t\Big )^p \D \sigma \\
&\, \lesssim \int_0^\infty \Big ( \sigma^{-(s+\alpha-(1-\mu+1/p))}   \int_0^\sigma t^{-1} \Big (t^{-(2-\mu)} \int_0^t \tau^{1-\mu} |u(t)-u(\tau)|_{D(A^s)} \D \tau\Big ) \D t\Big )^p \sigma^{-1}\D \sigma\\
&\,\lesssim \int_0^\infty  \sigma^{-p(s+\alpha-(1-\mu+1/p))}  \sigma^{-p(2-\mu)} \Big (\int_0^\sigma  \tau^{1-\mu} |u(\sigma)-u(\tau)|_{D(A^s)} \D \tau \Big )^p \sigma^{-1}\D \sigma\\
&\,\leq \int_0^\infty \int_0^\sigma  \tau^{p(1-\mu)} |u(\sigma)-u(\tau)|_{D(A^s)}^p \sigma^{-(1+p(s+\alpha))} \D \tau \D \sigma \leq [u]_{W_{p,\mu}^{s+\alpha}(\R_+;D(A^s))}^p.
\end{align*}
Thus (\ref{lg78}) holds. \eprf
 
From the above lemma we deduce a general trace theorem in the time variable
for the weighted anisotropic spaces. 

\begin{thm}\label{mott}\textsl{Let $E$ be a Banach space of class $\HT$, let $J=(0,T)$ be finite 
or infinite, and let $\Omega\subset \R^n$ be a bounded domain with smooth boundary, 
or $\Omega\in \{\R^n,\R_+^n\}$. Assume that $r\geq 0$, $\beta>0$, $0<1/p<\mu\le 1$ 
and  that $k\in \N_0$, $s\geq 0$ and $\alpha\in (0,2)$ satisfy
$$ s< k + 1-\mu+1/p < s+\alpha.$$  
Set  $H_{p,\mu}^s(W_p^r) := H_{p,\mu}^s(J; W_p^r(\Omega;E)),$ and analogously for the other 
anisotropic spaces. The order of differentiability of each
occurring $W_{p,\mu}$- and $W_p$-spaces is required to be noninteger. Then each of the spaces 
\beq\label{lg75}
H_{p,\mu}^{s+\alpha}(W_p^r) \cap H_{p,\mu}^s(W_p^{r+\beta}),\quad 
W_{p,\mu}^{s+\alpha}(H_p^r) \cap W_{p,\mu}^s(H_p^{r+\beta}), \quad 
W_{p,\mu}^{s+\alpha}(H_p^r) \cap H_{p,\mu}^s(W_p^{r+\beta}),
\eeq
is continuously embedded into
\beq\label{ts}
BU\!C^k\big (\ol{J}, B_{p,p}^{r+\beta(1+(s-(k+1-\mu+1/p))/\alpha)}(\Omega;E)\big ).
\eeq 
Moreover, for $\alpha \leq 1$ it holds 
\begin{align}\label{lg73}
W_{p,\mu}^\alpha(W_p^r)\cap L_{p,\mu}(W_p^{r+\beta})
  &\hra BU\!C\big(\ol{J}, B_{p,p}^{r+\beta(1-(1-\mu+1/p)/\alpha)}(\Omega;E)\big),\\
\label{lg74}
W_{p,\mu}^1(W_p^r)\cap W_{p,\mu}^s(W_p^{r+\beta})
  &\hra BU\!C\big(\ol{J}, B_{p,p}^{r+\beta(\mu-1/p)/(1-s)}(\Omega;E)\big).
\end{align}
All these embeddings remain true if one replaces $\Omega$ by its boundary $\partial\Omega$
or the $H_{p,\mu}$ and  $W_{p,\mu}$-spaces by ${}_0H_{p,\mu}$ and  ${}_0W_{p,\mu}$-spaces,
respectively. In the latter case, the embedding constants are independent of $T$ if
 $s+\alpha \leq 2$.}
\end{thm}
\bprf 
\textbf{(I)} 
We first indicate how to reduce the assertions to some basic cases.
Again we may assume that $J\times \Omega= \R_+\times \R^n$. The left translations
form a contractive $C_0$-semigroup on the five spaces under consideration, due to 
Lemma~\ref{sec:timeder}. Hence, one can see as in Proposition~III.1.4.2 of \cite{Ama95}
 that the asserted embeddings hold if we can show that the temporal trace operator 
 $\text{tr}_0 u = u|_{t=0}$ maps these function spaces continuously into 
 $$Y:=B_{p,p}^{r+\beta(1+(s-(k+1-\mu+1/p))/\alpha)}(\R^n;E),$$ 
 where one has to set $s=k=0$ for (\ref{lg73}) as well as $k=0$ and $\alpha = 1-s$ for (\ref{lg74}).
We may assume that 
$$s> k -\mu+1/p.$$ 
In fact, if this not the case, we take $\sigma\in(0,1)$
such that $s+\sigma\alpha> k -\mu+1/p$. Proposition~\ref{embani} then allows to embed,
say,  $Z:=W_{p,\mu}^{s+\alpha}(H_p^r) \cap 
H_{p,\mu}^s(W_p^{r+\beta})$ into $W^{s+\sigma\alpha}_{p,\mu}(H_p^{r+(1-\sigma)\beta})$.
For the space   $Z_1:= W_{p,\mu}^{s+\alpha}(H_p^r) 
\cap  W_{p,\mu}^{s+\sigma\alpha}(H_p^{r+(1-\sigma)\beta})$ we then have 
the embedding into \eqref{ts} proved below, which is precisely the assertion for $Z$.

Moreover, we only have to consider the case $k=0$ for the spaces in \eqref{lg75}. In fact, if $s\ge k$ we can 
simply apply $\partial_t^k$ to these spaces and then use the case $k=0$. If $s\in (k-\alpha,k)$,
then  we first use Proposition~\ref{embani} and embed, say,  $Z=W_{p,\mu}^{s+\alpha}(H_p^r) \cap 
H_{p,\mu}^s(W_p^{r+\beta})$ into $W^{k+\eps}_{p,\mu}(H_p^{r+(1-\sigma)\beta})$ with 
$\sigma=(k+\eps-s)/\alpha$, for some $\eps\in(0,s+\alpha-k)$. Now, $\partial_t^k$ maps
$Z$ into  $W_{p,\mu}^{s+\alpha-k}(H_p^r) \cap W_{p,\mu}^\eps(H_p^{r+(1-\sigma)\beta})$
and the trace result for $k=0$ will imply the assertion. Similarly, if we have $s+\alpha\in[1, 2-\mu+1/p]$, 
we can embed, say, $Z$  into  $Z_2:= W_{p,\mu}^{s+\alpha}(H_p^r) \cap 
W_{p,\mu}^{s+\gamma}(H_p^{r+(1-\sigma)\beta})$, where $\sigma=\gamma/\alpha$ and $s+\gamma\in (s,1)$.
The trace result for $Z_2$ will again imply the assertion for $Z$.

Finally, Proposition \ref{embani} yields the embeddings
\begin{align*}
H_{p,\mu}^{s+\alpha}(W_p^r) \cap H_{p,\mu}^s(W_p^{r+\beta}) &\hra 
W_{p,\mu}^{s+(1-\eps)\alpha}(H_p^{r+\eps\beta}) \cap W_{p,\mu}^{s+\eps \alpha}(H_p^{r+(1-\eps)\beta}),\\
W_{p,\mu}^{s+\alpha}(H_p^r) \cap H_{p,\mu}^s(W_p^{r+\beta}) &\hra 
W_{p,\mu}^{s+(1-\eps)\alpha}(H_p^{r+\eps\beta}) \cap W_{p,\mu}^{s+\eps \alpha}(H_p^{r+(1-\eps)\beta})
\end{align*}
where  $\eps>0$ is chosen such that the exponents of the $W$-spaces are noninteger. 
If we can show the assertion for the traces on the right hand side, then we also obtain the asserted
embeddings for the spaces on left hand side. So besides \eqref{lg73} and  \eqref{lg74}, it remains to 
consider the space 
\beq\label{WHcase}
W_{p,\mu}^{s+\alpha}(H_p^r) \cap W_{p,\mu}^s(H_p^{r+\beta})
 \qquad  \text{for \ $k=0$ \ and \ $s+\alpha<1$.}
\eeq

\textbf{(II)} We first show the trace result for the case \eqref{WHcase}. To that purpose,
we apply (\ref{lg78}) with 
$$X = H_p^{r-s\beta/\alpha}, \qquad A = (1-\Delta_n)^{\beta/2\alpha}, \qquad D(A) = H_p^{r+ (1-s)\beta/\alpha},$$ 
such that $D(A^s)  = H_p^r$ and $D(A^{s+\alpha}) = H_p^{r+\beta}.$ (Observe that $r-s\beta/\alpha<0$
is possible.) Hence $\tr_0$ is continuous
$$W_{p,\mu}^{s+\alpha}(H_p^r) \cap W_{p,\mu}^s(H_p^{r+\beta})\ra D_A\big (2s+\alpha- (1-\mu+1/p),p\big) 
= B_{p,p}^{r+\beta(1+(s-(1-\mu+1/p))/\alpha)},$$
as asserted. To obtain (\ref{lg73}), we use (\ref{lg76}) applied to
$$X = B_{p,p}^r, \qquad A = (1-\Delta_n)^{\beta/2\alpha}, \qquad D(A) = B_{p,p}^{r+\beta/\alpha},$$ 
giving $D_A(\alpha-(1-\mu+1/p),p) = B_{p,p}^{r+\beta(1-(1-\mu+1/p)/\alpha)}$.
Finally, we deduce (\ref{lg74}) from (\ref{lg77}) with
$$X = B_{p,p}^{r-s\beta/(1-s)}, \qquad A = (1-\Delta_n)^{\beta/2(1-s)}, \qquad D(A) = W_p^{r+\beta},$$
so that $D_A(1+s- (1-\mu+1/p),p) = B_{p,p}^{r+\beta(\mu-1/p)/(1-s)}$ in this case. \eprf


Arguing as in the proof of Theorem~4.5 in \cite{DSS08}, one should be able  to show that the 
temporal trace is surjective for all of the spaces under consideration in the Theorem \ref{mott}. 
At this point we only consider a right-inverse for  an important special case related to 
maximal $L_{p,\mu}$-regularity.

\begin{lem}\label{kv40}\textsl{Let $E$ be a Banach space of class $\HT$, $p\in (1,\infty)$, 
$\mu\in (1/p,1]$, $m\in \N$,  and let $\Omega\subset \R^n$ be a domain with compact smooth boundary 
$\partial\Omega$, or $\Omega\in \{\R^n,\R_+^n\}$. Then there exists a continuous right-inverse 
$$S: B_{p,p}^{2m(\mu-1/p)}(\Omega;E) \ra W_{p,\mu}^1\big (\R_+; L_p(\Omega;E)\big)
         \cap L_{p,\mu}\big(\R_+; W_p^{2m}(\Omega;E)\big)$$
of $\emph{\tr}_0$ which is given by
$$Su_0(t)= \calR_{\Omega} e^{-t (1-\Delta_n)^m} \calE_{\Omega}u_0, \qquad t>0.$$
Here, $\calE_{\Omega}$ is the extension to $\R^n$ from (\ref{ext1}) and $\calR_{\Omega}$ denotes 
the restriction from $\R^n$ to $\Omega$.}
\end{lem}

\bprf 
The operator  $A= (1-\Delta_n)^m$ with domain $W^{2m}_p(\R^n;E)$ generates
 an exponentially stable analytic $C_0$-semigroup  on $X=L_p(\R^n;E)$.  As in 
the proof of Lemma~\ref{wm22}, we obtain $D(A^\beta)= H_p^{2m\beta}(\R^n;E)$ for $\beta\in(0,1)$
so that $D_A(\mu-1/p,p)= B_{p,p}^{2m(\mu-1/p)}(\R^n;E)$ by reiteration. Theorem~1.14.5 of \cite{Tri94}
further implies that 
$$|e^{-\cdot A}x|_{W_{p,\mu}^{1}(\R_+;X)\cap L_{p,\mu}(\R_+; D(A))} \lesssim 
 |x|_{D_A(\mu-1/p,p)} \qquad  \text{for all \ } x\in D_A(\mu-1/p,p).$$ 
Using also  (\ref{ext1}), the assertion now follows.\eprf

We now specialize to weighted anisotropic spaces of the form 
\beq\label{wm17}
H_{p,\mu}^{s,2ms}(J\times \Omega;E):=H_{p,\mu}^s\big (J; L_p(\Omega;E)\big )
        \cap L_{p,\mu}\big (J;H_p^{2ms}(\Omega;E)\big ),
\eeq 
where $m\in \N$ and $s\in (0,1]$ satisfy $2ms\in \N$, and to the corresponding spaces
 where $H$ is replaced by $W$.  We are interested in the behavior of the outer normal derivative and 
more general boundary operators on the maximal $L_{p,\mu}$-regularity class. 
In view of Lemma \ref{spatiald}, we thus have to investigate the properties of the spatial trace 
$\text{tr}_\Omega$, i.e., 
$$\text{tr}_\Omega u := u|_{\partial\Omega},$$ 
on  spaces like (\ref{wm17}). We first give a heuristic derivation of the range space of
$\text{tr}_\Omega$.

It is known that $\text{tr}_\Omega$, which is originally only defined on $C_c^\infty(\R^n;E)$, extends 
uniquely to a continuous map 
\beq\label{kd200}
H_p^{2ms}(\Omega;E) \ra  W_{p}^{2ms-1/p}(\partial \Omega;E).
\eeq
This can be seen as in Theorems~2.9.3 and 4.7.1 of \cite{Tri94} for the scalar-valued case.  
Applied pointwise almost everywhere in time, $\text{tr}_\Omega$ then extends to a continuous map 
$$L_{p,\mu}\big(J; H_p^{2ms}(\Omega;E)\big) \ra  L_{p,\mu}\big(J; W_p^{2ms-1/p}(\partial \Omega;E)\big).$$ 
Observe that Proposition \ref{embani} yields  $H_{p,\mu}^{s,2ms}(J\times \Omega;E)
\hra H_{p,\mu}^{s-1/2mp}\big(J; H_p^{1/p}(\Omega;E)\big).$
Although $\tr_\Omega$ is not bounded from $H_p^{1/p}(\Omega;E)$ to $L_p(\partial\Omega;E)$, 
this  suggests that $\text{tr}_\Omega$ maps $H_{p,\mu}^{s,2ms}(J\times \Omega;E)$ also in a space like
$H_{p,\mu}^{s-1/2mp}\big(J; L_p(\partial\Omega;E)\big).$


To give a rigorous proof of the correct result, let us first assume that
$$J\times \Omega = \R_+\times \R_+^n,$$ 
and describe an alternative representation of $\tr_{\R_+^n}$. In the sequel we write  
$$x = (x',y)\in \R_+^n, \qquad x'\in \R^{n-1},\quad y\in \R_+.$$ 
Considering a function $u=u(t,x',y)$ on $\R_+\times \R_+^n$ as a function of $y\in \R_+$ with 
values in the functions of $(t,x')\in \R_+\times \R^{n-1}$, Fubini's theorem yields the embedding
$$\iota :L_{p,\mu}\big (\R_+; H_p^{2ms}(\R_+^n;E)\big)
      \hra  H_p^{2ms}\big(\R_+;L_{p,\mu}(\R_+;L_{p}(\R^{n-1};E))\big).$$ 
Since $2ms \geq 1$, the trace $\text{tr}_{0}:=\text{tr}_{y=0}$ thus acts on  
$L_{p,\mu}\big (\R_+; H_p^{2ms}(\R_+^n;E)\big)$ via $\text{tr}_{0}\circ \iota$ and maps this space 
continuously into $L_{p,\mu}\big(\R_+; L_{p}(\R^{n-1};E)\big)$. In a straightforward way, one can check 
that the set  $\trm{Step}\big(\R_+; C_c^\infty(\R^n;E)\big)$, consisting of step functions of the form 
$$
 \phi = \sum_{i=1}^l \alpha_i(\cdot) \phi_i, \qquad \alpha_i \in C_c(\R_+),\quad 
        \phi_i \in C_c^\infty(\R^n;E), \quad l\in \N,
$$
is dense in $L_{p,\mu}\big (J; H_p^{2m}(\R_+^n;E)\big)$, see Lemma~1.3.10 of \cite{Mey10}.
For $\phi \in \trm{Step}\big(\R_+; C_c^\infty(\R^n;E)\big)$ it clearly holds $\text{tr}_{\R_+^n}\phi 
= (\text{tr}_{0}\circ \iota )\phi$. We therefore obtain that
\beq\label{wm15}
\text{tr}_{\R_+^n} = \text{tr}_{0}\circ \iota \qquad \trm{ on } L_{p,\mu}\big(\R_+; H_p^{2ms}(\R_+^n;E)\big).
\eeq 
This representation allows to prove the temporal regularity for spatial traces of functions in 
$H_{p,\mu}^{s,2ms}\big (\R_+\times \R_+^n;E\big)$ as suggested above. Our proof is inspired by Lemma~3.5 of \cite{DHP07}.

\begin{lem} \label{sec:stt} \textsl{Let $E$ be a Banach space of class $\HT$,  $0<1/p<\mu\le1$, and let
$m\in\N$ and $s\in (0,1]$ satisfy $2ms\in \N$. Then the trace $\emph{\text{tr}}_{\R_+^n}$ maps continuously 
$$H_{p,\mu}^{s,2ms}(\R_+\times \R_+^n;E)\ra   W_{p,\mu}^{s-1/(2mp), 2ms -1/p}(\R_+\times \R^{n-1};E).$$ 
It is further surjective and has a continuous right-inverse.}
\end{lem}
\bprf Throughout this proof we set $\X:=L_{p,\mu}\big(\R_+; L_p(\R^{n-1};E)\big).$
 
\textbf{(I)} Considering a function in $L_{p,\mu}\big(\R_+; H_p^{2ms}(\R_+^n;E)\big)$ as a function of $y\in \R_+$ taking values in the functions of $(t,x')\in \R_+\times \R^{n-1}$, we obtain from Fubini's theorem that 
$$ L_{p,\mu}\big(\R_+; H_p^{2ms}(\R_+^n;E)\big) \hra H_p^{2ms}(\R_+;\X).$$ 
Moreover, it follows from $H_p^{2ms}(\R_+^n;E) \hra L_p\big(\R_+; H_p^{2ms}(\R^{n-1};E)\big)$  that 
$$L_{p,\mu}\big(\R_+; H_p^{2ms}(\R_+^n;E)\big) \hra L_p \big(\R_+; L_{p,\mu}(\R_+; H_p^{2ms}(\R^{n-1};E))\big).$$ 
Fubini's theorem and interpolation further yield
$$H_{p,\mu}^s\big(\R_+; L_p(\R_+^n;E)\big) \hra L_p\big(\R_+; H_{p,\mu}^s(\R_+; L_p(\R^{n-1};E))\big).$$ 
By Lemma \ref{wm22},  the operator $L= 1 - \partial_t + (-\Delta_{n-1})^m$ 
with domain $D(L)=H_{p,\mu}^{1,2m}(\R_+\times \R_+^n;E)$
is invertible on $\X$  and admits bounded imaginary powers with power angle not exceeding $\pi/2$. 
Hence for $\tau\in (0,1]$ its power $L^\tau$ has bounded imaginary powers with angle not larger than $\tau \pi/2$, and it holds 
\beq\label{kv34}
D(L^\tau) = H_{p,\mu}^\tau\big(\R_+; L_p(\R^{n-1};E)\big) \cap L_{p,\mu}\big(\R_+; H_p^{2m\tau}(\R^{n-1};E)\big),
\eeq 
again by Lemma \ref{wm22}. Therefore we obtain the embedding 
$$\wt{\iota}:H_{p,\mu}^{s,2ms}(\R_+\times \R_+^n;E) \hra H_p^{2ms}(\R_+; \X)\,\cap\, L_p(\R_+;D(L^s)),$$ 
and equation (\ref{wm15}) implies that $\text{tr}_{\R_+^n} = \text{tr}_{0} \circ \wt{\iota}$.

\textbf{(II)} We now claim that the space $H_p^{2ms}(\R; \X)\,\cap\, L_p(\R;D(L^s))$ embeds continuously into $H_p^1(\R; D(L^{s-1/2m}))$. To see this, we consider  the realization of the operators 
$A= (1-\partial_y^2)^{sm}$ and $B=L^s$ on $L_p(\R;\X)$ with domains 
$$D(A) = H_p^{2ms}(\R;\X) \qquad \text{and}\qquad D(B) =  L_p(\R; D(L^s)),$$ 
respectively. These operators are invertible and possess bounded imaginary powers with power angles 
equal to zero and $s\pi/2$, respectively, see e.g. Lemma \ref{wm22}.  Moreover, $A$ and $B$ are resolvent commuting on step functions in 
$L_p(\R;\X)$, which carries over to $L_p(\R;\X)$ by density. Thus Proposition \ref{dv} shows that the operator $A+B$ is invertible on $L_p(\R;\X)$ with domain 
$$D(A+B) = H_p^{2ms}(\R; \X)\,\cap\, L_p(\R;D(L^s)).$$ 
Since $A$, $B$ and $A+B$ are invertible, $|A^{1/2ms} B^{1-1/2ms} \cdot |_{L_p(\R;\X)}$ and 
$|(A+B) \cdot |_{L_p(\R;\X)}$ 
are equivalent norms on $D(A^{1/2ms} B^{1-1/2ms})= H_p^1\big (\R; D(L^{s-1/2m})\big )$ and $D(A+B)$, respectively. 
Now, Proposition~\ref{dv}  implies the asserted embedding.

\textbf{(III)} It follows from  extension and restriction that also 
$$H_p^{2ms}(\R_+; \X)\,\cap\, L_p(\R_+;D(L^s)) \hra H_p^1(\R_+; D(L^{s-1/2m})),$$ which implies that  $L^{s-1/2m}$ maps continuously $$H_p^{2ms}(\R_+; \X)\,\cap\, L_p(\R_+;D(L^s)) \ra H_p^{1}(\R_+; \X) \,\cap \, L_p(\R_+; D(L^{1/2m})).$$ 
Note that $L^{1/2m}$ is sectorial of angle at most $\pi/4m< \pi/2$, and thus $-L^{1/2m}$ generates an exponentially stable analytic $C_0$-semigroup on $\X$. Due to Theorem~III.4.10.2 in \cite{Ama95} we have
$$H_p^{1}(\R_+; \X) \,\cap\,L_p(\R_+; D(L^{1/2m})) \hra BU\!C\big([0,\infty); D_{L^{1/2m}}(1-1/p,p)\big),$$ 
and from the reiteration theorem we infer 
$$D_{L^{1/2m}}(1-1/p,p)  = D_{L}\big ((1-1/p)/2m,p\big ).$$ 

\textbf{(IV)} We now write 
$$\text{tr}_{\R_+^n}= \text{tr}_{0} \, L^{-(s-1/2m)} \, L^{s-1/2m} \, \wt{\iota},$$
where $L^{s-1/2m}$ and its inverse are applied pointwise.
By the above considerations, the operator $L^{s-1/2m} \wt{\iota}$ is continuous 
$$H_{p,\mu}^{s,2ms}(\R_+\times \R_+^n;E) \ra BU\!C\big([0,\infty); D_{L}((1-1/p)/2m,p)\big).$$ 
Clearly, $\text{tr}_{0}$ and $L^{-(s-1/2m)}$ commute on $BU\!C\big([0,\infty); D_{L}((1-1/p)/2m,p)\big)$.
By reiteration $L^{-(s-1/2m)}$ maps $D_{L}((1-1/p)/2m,p)$ continuously into
 $$D_{L}(s-1/2mp,p) = W_{p,\mu}^{s-1/(2mp), 2ms -1/p}(\R_+\times \R^{n-1};E),$$
where we use Lemma~\ref{wm22}. Hence, $\text{tr}_{\R_+^n}$ maps as asserted. 

\textbf{(V)} A continuous right-inverse of $\tr_{\R_+^n}$ is defined by
$$e^{-y L^{1/2m}} g(t,x'), \qquad t>0, \quad (x',y)\in \R_+^n.$$ 
for $g\in W_{p,\mu}^{s-1/(2mp), 2ms -1/p}(\R_+\times \R^{n-1};E)$.
To see this, observe that 
$$
W_{p,\mu}^{s-1/(2mp), 2ms -1/p}(\R_+\times \R^{n-1};E) = D_{L^{1/2m}}(2ms -1/p,p)
$$
by reiteration. For $f\in  D_{L^{1/2m}}(2ms -1/p,p)$ and $k=2ms-1$, we have that 
$(-\partial_y)^ke^{-y L^{1/2m}}f= e^{-y L^{1/2m}} L^{k/2m}f$ and $ L^{k/2m}f\in D_{L^{1/2m}}(1 -1/p,p)$.
Using e.g.\ Theorem~1.14.5 of \cite{Tri94}, we then derive that
 $e^{-\cdot L^{1/2m}}$ maps $D_{L^{1/2m}}(2ms -1/p,p)$ continuously into 
$\Y := H_{p}^{2ms}(\R_+; \X) \, \cap \, L_p(\R_+; D(L^s)).$
Therefore (\ref{kv34}) yields
$$\Y \hra L_p\big (\R_+; H_{p,\mu}^s(\R_+; L_p(\R^{n-1};E))\big) \hra H_{p,\mu}^s\big(\R_+; L_p(\R_+^{n};E)\big).$$ 
Moreover,
\begin{align*}
\Y&\, \hra  H_{p}^{2ms}\big(\R_+; L_{p,\mu}(\R_+;L_p(\R^{n-1};E))\big) \, \cap \, L_p\big(\R_+, L_{p,\mu}(\R_+, H_p^{2ms}(\R^{n-1};E))\big)\\
&\,\hra L_{p,\mu}\big(\R_+; H_p^{2ms}(\R_+^n;E)\big),
\end{align*}
where the latter embedding follows again from Fubini's theorem. This shows the asserted continuity of the right-inverse. \eprf

The mapping properties of the spatial trace on a domain are now obtained from the half-space case via localization.

\begin{thm}\label{sec:generalstt}\textsl{Let $E$ be a Banach space of class $\HT$, $J=(0,T)$ be a finite 
or infinite, $0<1/p<\mu\le1$, and let $m\in\N$ and $s\in (0,1]$ satisfy $2ms\in \N$. 
Assume that  $\Omega\subset \R^n$ is a domain with compact smooth boundary, or $\Omega\in \{\R^n,\R_+^n\}$. 
Then  the trace $\emph{\text{tr}}_{\Omega}$ is continuous 
$$H_{p,\mu}^{s,2ms}(J\times \Omega;E) \ra   W_{p,\mu}^{s-1/2mp, 2ms -1/p}(J\times \partial \Omega;E).$$ 
It is further surjective and has a continuous right-inverse. The operator norm of $\emph{\text{tr}}_{\Omega}$ on ${}_0H_{p,\mu}^{s,2ms}(J\times \Omega;E)$ is independent of the length of $J$.}
\end{thm}
\bprf Using the extension operators $\calE_J$ and $\calE_J^0$ from Lemma \ref{sec:ext3}, it suffices to consider the case $J=\R_+$. We describe $\partial\Omega$ by a finite number of charts $(U_i, \vphi_i)$ and a partition of unity $\{\psi_i\}$ subordinate to the cover $\bigcup_i U_i$. We further denote by $\Phi_i$ the push-forward with respect to $\vphi_i$; i.e., $\Phi_i u = u\circ \vphi_i^{-1}$. Then for $\phi\in \trm{Step}\big (\R_+; C_c^\infty(\R^n;E)\big )$ it holds 
\beq\label{lg70}
\text{tr}_{\Omega}\phi =  \sum_i \nolimits \Phi_i^{-1}\big (\text{tr}_{\R_+^n} \Phi_i (\psi_i\phi)\big ) \qquad \trm{on } \partial \Omega.
\eeq 
By restriction to $\Omega\cap U_i$ and trivial extension from $\R_+^n \cap \vphi_i(U_i)$ to $\R_+^n$, for each $i$ we obtain that   $\Phi_i (\psi_i\cdot )$ maps continuously 
$$H_{p,\mu}^{s,2ms}(\R_+\times \Omega;E)\ra H_{p,\mu}^{s,2ms}(\R_+\times \R_+^n;E).$$ 
Applying Lemma \ref{sec:stt} and restricting back to $\R_+^n \cap \vphi_i(U_i)$ yields that  
$\Phi_i^{-1} \text{tr}_{\R_+^n}$ maps the latter space continuously into 
$W_{p,\mu}^{s-1/2mp, 2ms -1/p}(\R_+\times \partial \Omega ;E).$
Thus  the map $\sum_i \Phi_i^{-1}\big (\text{tr}_{\R_+^n} \Phi_i (\psi_i \cdot)\big )$ is continuous
$$H_{p,\mu}^{s,2ms}(\R_+\times \Omega;E) \ra   W_{p,\mu}^{s-1/2mp, 2ms -1/p}(\R_+\times \partial \Omega;E).$$ 
 Since $\trm{Step}\big (\R_+; C_c^\infty(\R^n;E)\big )$ is dense in $L_{p,\mu}\big (\R_+; H_p^{2ms}(\Omega;E)\big)$, the representation (\ref{lg70}) holds for all elements of this space, and in particular for all 
functions from $H_{p,\mu}^{s,2ms}(J\times \Omega;E)$. This shows that $\tr_\Omega$ is continuous as asserted. Similiar localization arguments also reduce the existence of a  continuous right-inverse 
to the half-space case treated in Lemma~\ref{sec:stt}. \eprf

\end{document}